\documentclass[12pt,english]{article}
\usepackage[T1]{fontenc}
\usepackage[latin9]{inputenc}
\usepackage[letterpaper]{geometry}
\usepackage{mathrsfs}
\usepackage{amsmath}
\usepackage{amsthm}
\usepackage{amssymb}
\usepackage{graphicx}

\makeatletter
\numberwithin{equation}{section}
\numberwithin{figure}{section}
\theoremstyle{definition}
\newtheorem*{defn*}{\protect\definitionname}
\theoremstyle{plain}
\newtheorem*{thm*}{\protect\theoremname}
\theoremstyle{plain}
\newtheorem*{cor*}{\protect\corollaryname}
\theoremstyle{plain}
\newtheorem*{lem*}{\protect\lemmaname}
\theoremstyle{definition}
\newtheorem*{example*}{\protect\examplename}

\@ifundefined{date}{}{\date{}}
\usepackage[bottom]{footmisc}

\makeatother

\usepackage{babel}
\providecommand{\corollaryname}{Corollary}
\providecommand{\definitionname}{Definition}
\providecommand{\examplename}{Example}
\providecommand{\lemmaname}{Lemma}
\providecommand{\theoremname}{Theorem}

\begin{document}
\title{\noindent \textbf{Constructive Projective Geometry}}
\author{Mark Mandelkern}

\maketitle
\vspace{1cm}

\tableofcontents{}

\vspace{1cm}

\section{Introduction \label{Introduction}}

\noindent Of the great theories of classical mathematics, projective
geometry, with its powerful concepts of symmetry and duality, has
been exceptional in continuing to intrigue investigators. The challenge
put forth by Errett Bishop (1928-1983),
\begin{quotation}
Every theorem proved with {[}nonconstructive{]} methods presents a
challenge: to find a constructive version, and to give it a constructive
proof. {[}B67, p. x; BB85, p. 3{]}
\end{quotation}
\noindent and Bishop's ``Constructivist Manifesto'' {[}B67, Chapter
1; BB85, Chapter 1{]}, motivate a large portion of current constructive
work. This challenge can be answered by discovering the hidden constructive
content of classical projective geometry. Here we briefly outline,
with few details, recent constructive work on the real projective
plane, and  projective extensions of affine planes. Special note is
taken of a number of interesting open problems that remain; these
show that constructive projective geometry  is still a theory very
much in need of further effort. 

There has been a considerable amount of work in the constructivization
of geometry, on various topics, in different directions, and from
diverse standpoints. For the constructive extension of an affine plane
to a projective plane, see {[}H59, vDal63, M13, M14{]}. For the constructive
coördinatization of a plane, see {[}M07{]}. Intuitionistic axioms
for projective geometry were introduced by A. Heyting {[}H28{]}, with
further work by D. van Dalen {[}vDal96{]}. Work in constructive geometry
by M. Beeson {[}Be10, Be16{]} uses Markov's Principle, a nonconstructive
principle which is accepted in recursive function theory, but not
in the Bishop-type strict constructivism that is adhered to in the
present paper.  M. Lombard and R. Vesley {[}LV98{]} construct an axiom
system for intuitionistic plane geometry, and study it with the aid
of recursive function theory. The work of J. von Plato {[}vPla95,
vPla98, vPla10{]} in constructive geometry, proceeding from the viewpoint
of formal logic, is related to type theory, computer implementation,
and combinatorial analysis. The work of V. Pambuccian, e.g., {[}Pam98,
Pam01, Pam03, Pam05, Pam11{]}, also proceeding within formal logic,
covers a wide range of topics concerning axioms for constructive geometry. 

The Bishop-type constructive mathematics discussed in the present
paper proceeds from a viewpoint well-nigh opposite that of either
formal logic or recursive function theory. For further details concerning
this distinction, see {[}B65, B67, B73, B75, BB85, BR87{]}. 

\part{Real projective plane\label{PART I - real projective plane}}

Arend Heyting (1898-1980), in his doctoral dissertation {[}H28{]},
began the constructivization of projective geometry. Heyting's work
involves both synthetic and analytic theories. Axioms for projective
space are adopted; since a plane is then embedded in a space of higher
dimension, it is possible to include a proof of Desargues\textquoteright s
Theorem. For the coördinatization of projective space, axioms of order
and continuity are assumed. The theory of linear equations is included,
and results in analytic geometry are obtained. Later, Heyting discussed
the role of axiomatics in constructive mathematics as follows: 
\begin{quotation}
At first sight it may appear that the axiomatic method cannot be used
in intuitionistic mathematics, because there are only considered mathematical
objects which have been constructed, so that it makes no sense to
derive consequences from hypotheses which are not yet realized. Yet
the inspection of the methods which are actually used in intuitionistic
mathematics shows us that they are for an important part axiomatic
in nature, though the significance of the axiomatic method is perhaps
somewhat different from that which it has in classical mathematics.
{[}H59, p. 160{]} 
\end{quotation}
Recent work {[}M16, M18{]}, briefly outlined below, constructivized
the synthetic theory of the real projective plane as far as harmonic
conjugates, projectivities, the axis of homology, conics, Pascal's
Theorem, and polarity. Axioms only for a plane are used. The basis
for the constructivization is the extensive literature concerning
the classical theory, including works of O. Veblen and J. W. Young
{[}VY10, Y30{]}, H. S. M. Coxeter {[}Cox55{]}, D. N. Lehmer {[}Leh17{]},
L. Cremona\textsc{ }{[}Cre7\textsc{3{]}, }and G. Pickert {[}Pic75{]}.
An entertaining history of the classical theory is found in Lehmer's
last chapter{\small{}. }{\small\par}

\section{Axioms\label{SEC - Axioms}}

For nearly two hundred years a sporadic and sometimes bitter debate
has continued, concerning the value of synthetic versus analytic methods.
In his Erlangen program of 1872, Felix Klein sought to mediate the
dispute: 
\begin{quotation}
The distinction between modern synthetic and modern analytic geometry
must no longer be regarded as essential, inasmuch as both subject-matter
and methods of reasoning have gradually taken a similar form in both.
We chose therefore as common designation of them both the term \emph{projective
geometry.} Although the synthetic method has more to do with space-perception
and thereby imparts a rare charm to its first simple developments,
the realm of space-perception is nevertheless not closed to the analytic
method, and the formulae of analytic geometry can be looked upon as
a precise and perspicuous statement of geometrical relations. On the
other hand, the advantage to original research of a well formulated
analysis should not be underestimated, - an advantage due to its moving,
so to speak, in advance of the thought. But it should always be insisted
that a mathematical subject is not to be considered exhausted until
it has become intuitively evident, and the progress made by the aid
of analysis is only a first, though a very important, step. {[}Kle72{]} 
\end{quotation}
In the synthetic work summarized below, axioms are formulated which
can be traced to an analytic model based on constructive  properties
of the real numbers, and the resulting axiom system is used to construct
a synthetic projective plane $\mathbb{P}$. In this sense, the construction
of the plane $\mathbb{P}$ takes into account Bishop's thesis: ``All
mathematics should have numerical meaning'' {[}B67, p. ix; BB85,
p. 3{]}. 

\subsection{Axiom Group C\label{Axiom-Group-C}}

The constructivization of {[}M16{]}, resulting in the projective plane
$\mathbb{P}$, uses only axioms for a plane. There exist non-Desarguesian
projective planes; see, for example, {[}Wei07{]}. This means that
Desargues\textquoteright s Theorem must be taken as an axiom; it is
required to establish the essential properties of harmonic conjugates.
Other special features of the axiom system are also required, to obtain
constructive versions of the most important classical results. The
consistency of the axiom system is verified by means of an analytic
model, discussed below in Section \ref{SEC 8 Consistency}; the properties
of this model have guided the choice of axioms. 

The constructive axiom group C, adopted for the projective plane $\mathbb{P}$
in {[}M16, Section 2{]}, has seven initial axioms. The first four
are those usually seen for a classical projective plane; e.g., two
points determine a line, and two lines intersect at a point. The last
three axioms, which have special constructive significance, will be
discussed below. 

For the construction of the projective plane $\mathbb{P}$, there
is given a family $\mathscr{P}$ of points and a family $\mathscr{L}$
of lines, along with equality and inequality relations for each family.
The \emph{inequality relations }assumed for the families $\mathscr{P}$
and $\mathscr{L}$, both denoted $\neq$, are \emph{tight apartness
relation}s; thus, for any elements $x,y,z$, the following conditions
are satisfied: 

(i) $\neg(x\neq x)$. 

(ii) If $x\neq y$, then $y\neq x$. 

(iii) If $x\neq y$, then either $z\neq x$ or $z\neq y$. 

(iv) If $\neg(x\neq y)$, then $x=y$. 

\noindent The notion of an apartness relation was introduced by Brouwer
{[}Brou24{]}, and developed further by Heyting {[}H66{]}. Property
(iii) is known as \emph{cotransitivity}, and (iv) as \emph{tightness}.
The implication ``$\neg(x=y)$ implies $x\neq y$'' is invalid in
virtually all constructive theories, the inequality being the stronger
of the two conditions. For example, with real numbers considered constructively,
$x\neq0$ means that there exists an integer $n$ such that $1/n<|x|$,
whereas $x=0$ means merely that it is contradictory that such an
integer exists. For more details concerning the constructive properties
of the real numbers, see {[}B67, BB85, BV06{]}; for a comprehensive
treatment of constructive inequality relations, see {[}BR87, Section
1.2{]}. 

A given \emph{incidence relation,} written $P\in l$, links the two
families; we say that\emph{ the point $P$ lies on the line $l$},
or that $l$ \emph{passes through} $P$. A line is not viewed as a
set of points; the set $\overline{l}$ of points that lie on a line
$l$ is a \emph{range of points}, while the set $Q^{*}$ of lines
that pass through a point $Q$ is a \emph{pencil of lines}. The outside
relation $P\notin l$ is obtained by a definition: 
\begin{defn*}
\emph{\label{Defn Outside relation}}\textsc{Outside relation. }For
any point $P$\emph{ }on the projective plane \emph{$\mathbb{P}$,}
and any line $l$, it is said that \emph{$P$} \emph{lies outside
}$l$ (and\emph{ $l$ avoids $P$}), and written $P\notin l$, if
$P\neq Q$ for all points $Q$ that lie on $l$. {[}M16, Defn. 2.3{]} 
\end{defn*}
This condition for the relation $P\notin l$, when viewed classically,
is simply the negation of the condition $P\in l$, when written as
the tautology ``there exists $Q\in l$ such that $P=Q$''. Constructively,
however, the condition acquires a strong, positive significance, derived
from the character of the condition $P\neq Q$. 

Several axioms connect these relations: \\

\noindent \textbf{Axiom C5.} \emph{\label{Axiom C5}For any lines
$l$ and $m$ on the projective plane $\mathbb{P}$, if there exists
a point $P$ such that $P\in l$, and $P\notin m$, then $l\neq m$.
}\\

The implication ``If $\neg(P\in l)$, then $P\notin l$'' is nonconstructive.
However, we have: \\

\noindent \textbf{Axiom C6.} \emph{\label{Axiom C6}For any point
$P$ on the projective plane $\mathbb{P}$, and any line $l$, if
$\neg(P\notin l)$, then $P\in l$.}\\

Axiom C6 would be immediate in a classical setting, where\emph{ $P\notin l$}
means $\mbox{\ensuremath{\neg(P\in l)}}$; $\mbox{applying}$ the
law of excluded middle, a double negation results in an affirmative
statement. For a constructive treatment, where the condition \emph{$P\notin l$}
is not defined by negation, but rather by the affirmative definition
above, Axiom C6 must be assumed; it is analogous to the tightness
property of the inequality relations that are assumed for points and
lines.

For the metric plane $\mathbb{R}^{2}$, the condition of Axiom C6
follows from the analogous constructive property of the real numbers:
``For any real number $\alpha$, if $\neg(\alpha\ne0)$, then $\alpha=0$'',
interpreting the outside relation in terms of distance. For the analytic
model $\mathbb{P}^{2}(\mathbb{R})$, which motivates the axiom system,
Axiom C6 is verified using this constructive property of the real
numbers. 

The following axiom has a preëminent standing in the axiom system;
it is indispensable for virtually all constructive proofs involving
the  projective plane $\mathbb{P}$. The point of intersection of
distinct lines \emph{$l$ }and\emph{ $m$} is denoted \emph{$l\cdot m$}.
\\

\noindent \textbf{Axiom C7.}\label{Axiom C7} \emph{If $l$ and $m$
are distinct lines on the projective plane $\mathbb{P}$, and $P$
is a point such that $P\neq l\cdot m$, then either $P\notin l$ or
$P\notin m$. }\\

This axiom is a strongly worded, yet classically equivalent, constructive
form of a classical axiom: ``distinct lines have a \emph{unique}
common point'', which means only that if the points $P$ and $Q$
both lie on both lines, then $P=Q$. Axiom C7, a (classical) contrapositive
of the classical axiom, is significantly stronger, since the condition
$P\notin l$ is an affirmative condition. 

Heyting and van Dalen have used an apparently weaker version of Axiom
C7; it is Heyting's Axiom VI {[}H28{]}, and van Dalen's Lemma 3(f),
obtained using his axiom Ax5 {[}vDal96{]}.   This weaker version states:
``If $l$ and $m$ are distinct lines, $P$ is a point such that
$P\neq l\cdot m$, and $P\in l$, then $P\notin m$.'' However, it
is easily shown that the two versions are equivalent.

Axiom C4 states that at least three distinct points lie on any given
line; this is the usual classical axiom.\emph{ }Then, for the study
of projectivities, Axiom E is added, increasing the required number
of points to six. Recently, a constructive proof in {[}M18{]}, of
an essential result concerning harmonic conjugates, required at least
eight points on a line; thus we have: \\

\emph{Problem. }Determine the  minimum number of points on a line
that are required for the various constructive proofs concerning the
 projective plane $\mathbb{P}$. Examine the propositions involved
for the exceptional small finite planes. \\

The axioms and definitions of constructive projective geometry can
be given a variety of different arrangements. For example, in {[}vDal96{]}
the outside relation $P\notin l$ is taken as a primitive notion,
and the condition of Axiom C6 above becomes the definition of the
incidence relation $P\in l$. See also {[}H28, Pam05, Pam11, vPla95{]}. 

The axiom system could be extended; thus we have: \\

\emph{Problem. }Extend the constructive axiom group C to projective
space, and derive constructive versions of the main classical theorems. 

\subsection{Desargues's Theorem}

Desargues's Theorem is assumed as an axiom; the converse is then proved
as a consequence.

Two triangles are \emph{distinct} if corresponding vertices are distinct
and corresponding sides are distinct; it is then easily shown that
the three lines joining corresponding vertices are distinct, and the
three points of intersection of corresponding sides are distinct.
Distinct triangles are said to be \emph{perspective from the center
$O$} if the lines joining corresponding vertices are concurrent at
the point $O$, and $O$ lies outside each of the six sides. Distinct
triangles are said to be \emph{perspective from the axis $l$} if
the points of intersection of corresponding sides are collinear on
the line $l$, and $l$ avoids each of the six vertices. \\

\noindent \textbf{Axiom D.}\emph{\label{Axiom D} }\textsc{\small{}Desargues's
Theorem.}\emph{ If two triangles are perspective from a center, then
they are perspective from an axis. }\\

The proof of the converse is included below, as an example of constructive
methods in geometry. 

\includegraphics[viewport=50bp 20bp 950bp 390bp,clip,scale=0.48]{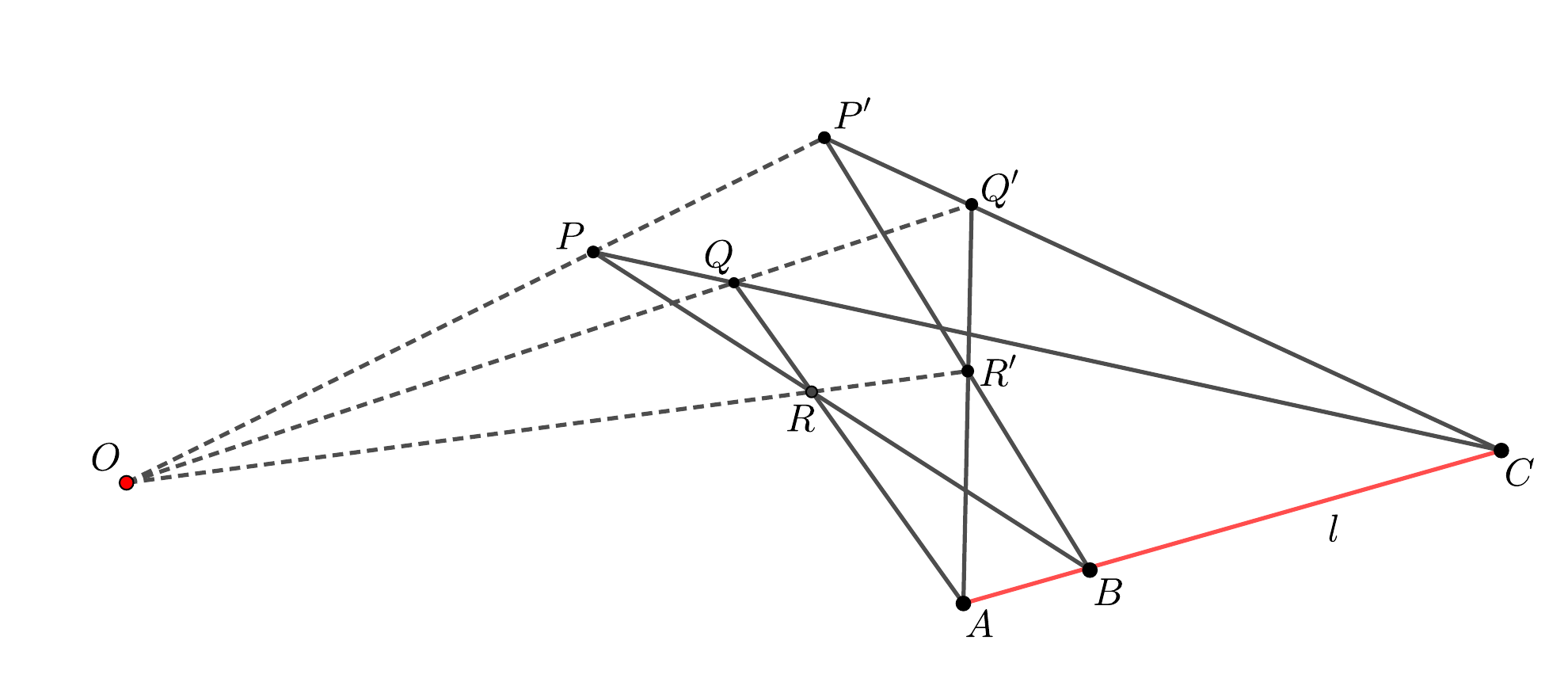}
\begin{thm*}
\label{Desargues converse}\textsc{\small{}Converse of Desargues's
Theorem.} If two triangles are perspective from an axis, then they
are perspective from a center.\emph{ {[}M16, Thm. 3.2{]}}
\end{thm*}
\begin{proof}
We are given distinct triangles $PQR$ and $P'Q'R'$, with collinear
points of intersection of corresponding sides, $A=QR\cdot Q'R'$,
$B=PR\cdot P'R'$, $C=PQ\cdot P'Q'$, and with all six vertices lying
outside the axis $l=AB$. Set $O=PP'\cdot QQ'$. 

The points $A,Q,Q'$ are distinct, and the points $B,P,P'$ are distinct.
Since $Q\neq A=QR\cdot Q'R'$, it follows from Axiom C7 that $Q\notin Q'R'=AQ'$;
thus the points $A,Q,Q'$ are noncollinear, and similarly for $B,P,P'$.
Since $P\notin AB$, we have $AB\neq BP$. Since $A\neq B=AB\cdot BP$,
it follows that $A\notin BP$, so $AQ\neq BP$. By symmetry, $AQ'\neq BP'$.
This shows that the triangles $AQQ'$, $BPP'$ are distinct. 

The lines $AB$, $PQ$, $P'Q'$, joining corresponding vertices of
the triangles $AQQ'$, $BPP'$, are concurrent at $C$. Since $Q\notin AB$,
we have $AB\ne AQ$. From $C\neq A=AB\cdot AQ$, it follows that $C\notin AQ$;
by symmetry, $C\notin AQ'$. Since $Q'\neq C=PQ\cdot P'Q'$, it follows
that $Q'\notin PQ$; thus $QQ'\neq PQ$, i.e., $CQ\neq QQ'$. From
$C\neq Q=CQ\cdot QQ'$, we have $C\notin QQ'$. Thus $C$ lies outside
each side of triangle $AQQ'$, and similarly for triangle $BPP'$.
Thus the triangles $AQQ'$, $BPP'$ are perspective from the center
$C$. 

It follows from Axiom D that the triangles $AQQ'$, $BPP'$ are perspective
from the axis $(AQ\cdot BP)(AQ'\cdot BP')=RR'$, the axis avoids all
six vertices, and $O\in RR'$. Thus the lines $PP'$, $QQ'$, $RR'$,
joining corresponding vertices of the given triangles, are concurrent
at $O$. Since $Q\notin RR'$, we have $Q\neq O$. From $O\neq Q=QQ'\cdot PQ$,
it follows that $O\notin PQ$. By symmetry, $O$ lies outside each
side of the given triangles. 

Hence the triangles $PQR$ and $P'Q'R'$ are perspective from the
center $O$. 
\end{proof}

\subsection{Duality}

Given any statement, the \emph{dual statement} is obtained by interchanging
the words ``point'' and ``line''. For example: \\

\noindent \textbf{Dual of Axiom C5}\textbf{\emph{. }}\emph{For any
points $P$ and $Q$ on the projective plane $\mathbb{P}$, if there
exists a line l such that $P\in l$, and $Q\notin l$, then $P\neq Q$.
}\\

\noindent \textbf{Dual of Axiom C7. \label{Axiom C7 Dual}}\emph{Let
$A$ and $B$ be distinct points on the projective plane $\mathbb{P}$.
If $l$ is a line such that $l\neq AB$, then either $A\notin l$
or $B\notin l$.} \\

\noindent Clearly, Axiom C6 is self-dual.\emph{ Duality} in a given
system is the principle that the dual of any true statement is also
true. Duality of the construction of the  plane $\mathbb{P}$, and
of the axiom system, is verified as follows: 
\begin{thm*}
The definition of the  projective plane $\mathbb{P}$ is self-dual.
The dual of each axiom in axiom group C is valid on $\mathbb{P}$.
\emph{{[}M16, Thm. 2.10{]}}
\end{thm*}
The dual of the definition of the outside relation $P\notin l$ is
also verified: 
\begin{thm*}
\label{Thm dual outside}Let $P$ be any point on the projective plane
$\mathbb{P}$, and $l$ any line. Then $P\notin l$ if and only if
$l\neq m$ for any line $m$ that passes through $P$. \emph{{[}M16,
Thm. 2.11{]}}
\end{thm*}

\section{Harmonic conjugates \label{SEC - Harmonic-conjugates}}

\noindent In the construction of the  projective plane $\mathbb{P}$,
harmonic conjugates have an essential role, with applications to projectivities,
involutions, and polarity. In the drawing below, the quadrangle $PQRS$,
which is often used classically, appears to determine the harmonic
conjugate $D$ of the point $C$, with respect to the base points
$A$ and $B$. However, this is only valid when $C$ is distinct from
each base point; thus we must use a definition that applies to every
point on the base line $AB$. 

\includegraphics[clip,scale=0.78]{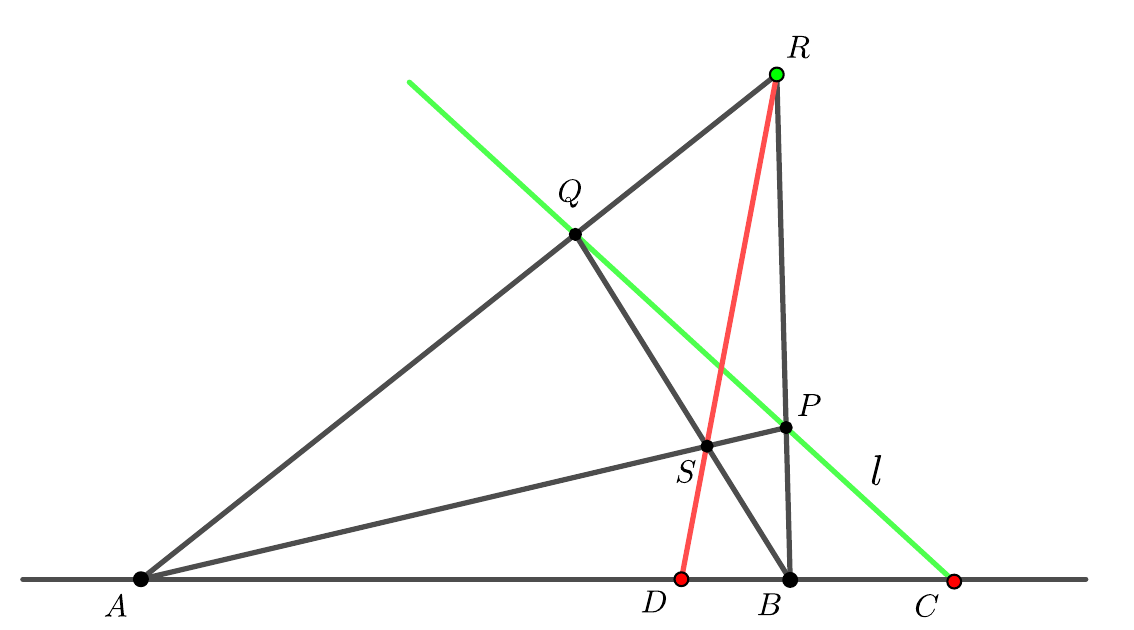}
\begin{defn*}
\label{Defn Harm Conj}Let\emph{ $A$} and \emph{$B$} be distinct
points on the projective plane $\mathbb{P}$. For any point \emph{$C$}
on the line $AB$, select a line \emph{$l$} through $C$, distinct
from $AB$, and select a point \emph{$R$} lying outside each of the
lines \emph{$AB$} and $l$. Set $P=BR\cdot l$, $Q=AR\cdot l$, and
\emph{$S=AP\cdot BQ$}. The point \emph{$D=AB\cdot RS$} is called
the \emph{harmonic conjugate} \emph{of} $C$ \emph{with respect to
the points} $A,B$; we write $D=h(A,B;C)$. {[}M16, Defn. 4.1{]}
\end{defn*}
Since the construction of a harmonic conjugate requires the selection
of auxiliary elements, it must be verified that the result is independent
of the choice of these auxiliary elements. The proof given in {[}M16{]}
for the invariance theorem is incorrect; apart from the error, the
proof there is excessively complicated, and objectionable on several
counts. A correct proof appears in a later paper.
\begin{thm*}
\textsc{\small{}Invariance Theorem}\textsc{\emph{\small{}.}} Let $C$
be any point on the line $AB$, and let auxiliary element selections
$(l,R)$ and $(l',R')$ be used to construct harmonic conjugates $D$
and $D'$ of the point $C$. Then $D=D'$; the harmonic conjugate
construction is independent of the choice of auxiliary  elements.
\emph{{[}M18, Thm. 3.2{]} }
\end{thm*}
In the special case of a point distinct from both base points, constructive
harmonic conjugates can be related to the traditional quadrangle configuration,
due to Philippe de La Hire (1640 \textendash{} 1718): 
\begin{cor*}
Let $A,B,C,D$ be collinear points, with $A\neq B$, and $C$ distinct
from each of the points $A$ and $B$. Then $D=h(A,B;C$) if and only
if there exists a quadrangle with vertices outside the line $AB$,
of which two opposite sides intersect at $A$, two other opposite
sides intersect at $B$, while the remaining two sides meet the base
line $AB$ at $C$ and $D$. \emph{{[}M18, Cor. 3.3{]} }
\end{cor*}

\section{Projectivities\label{SEC - Projectivities}}

The elementary mappings of a projective plane are \emph{sections,}
bijections relating a pencil of lines with a range of points. Certain
combinations of sections result in \emph{projections,} mapping a range
of points onto another range, projecting from a center, or mapping
a pencil of lines onto another pencil, projecting from an axis. These
sections and projections are the \emph{perspectivities} of the plane.

\includegraphics[viewport=-30bp 0bp 1675bp 499bp,clip,scale=0.25]{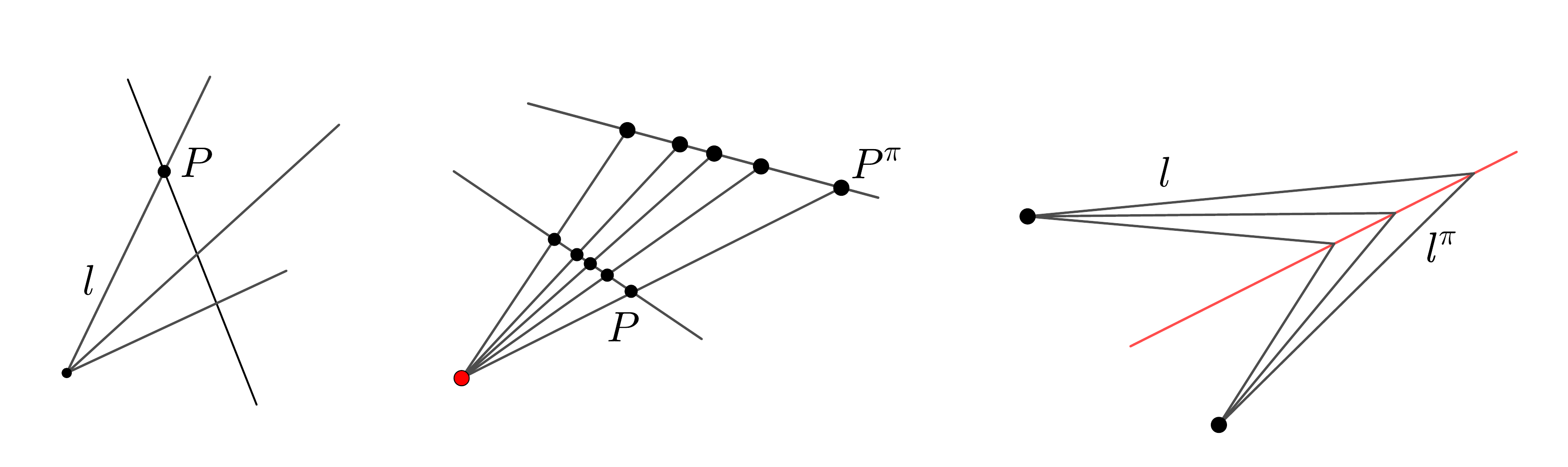}

The product (composition) of two perspectivities need not be a perspectivity.
For the  projective plane $\mathbb{P}$, a finite product of perspectivities
is called a \emph{projectivity;} this is the definition used by Jean-Victor
Poncelet (1788 \textendash{} 1867) {[}Pon22{]}. Subsequently, Karl
Georg Christian von Staudt (1798 \textendash{} 1867) {[}vSta47{]}
defined a\emph{ projectivity }as a mapping of a range or a pencil
that preserves harmonic conjugates. Classically, the two notions of
perspectivity are equivalent. Constructively, we have:
\begin{thm*}
A projectivity of the  projective plane $\mathbb{P}$ preserves harmonic
conjugates. Thus every Poncelet projectivity is a von Staudt projectivity.
\emph{{[}M16, Thm. 5.3{]}}
\end{thm*}
However, the constructive content of the converse is not known; thus
we have: \\

\emph{Problem.}\label{Prob Poncelet vonStaudt} On the  projective
plane $\mathbb{P}$, show that every von Staudt projectivity is a
Poncelet projectivity, or construct a counterexample. \\

It is necessary to establish the existence of projectivities:
\begin{thm*}
\label{Projectivity existence} Given any three distinct points $P,Q,R$
in a range $\overline{l}$, and any three distinct points $P',Q',R'$
in a range $\overline{m}$, there exists a projectivity $\pi:\overline{l}\rightarrow\overline{m}$
such that the points $P$,$Q$,$R$ map into the points $P',Q',R'$,
in the order given. \emph{{[}M16, Thm. 5.6{]}}
\end{thm*}
Classically, the projectivity produced by this theorem is the product
of at most three perspectivities. However, the constructive proof
in {[}M16{]} requires six perspectivities; thus we have: \\

\emph{Problem. }Determine the minimum number of perspectivities required
for the above theorem. \\

A projectivity $\pi$ of order 2 ($\pi^{2}$ is the identity $\iota$)
is called an \emph{involution}; this term was first used by Girard
Desargues (1591-1661). In {[}Des64{]}, Desargues introduced seventy
new geometric terms; they were considered highly unusual, and met
with sharp criticism and ridicule by his contemporaries. Of these
seventy terms, \emph{involution} is the only one to have survived.
 One example of an involution is the harmonic conjugate relation: 
\begin{thm*}
Let $A$ and $B$ be distinct points in a range $\overline{l}$, and
let $\upsilon$ be the mapping of harmonic conjugacy with respect
to the base points $A,B$; i.e., set $X^{\upsilon}=h(A,B;X)$, for
all points $X$ in the range $\overline{l}$. Then $\upsilon$ is
an involution. \emph{{[}M16, Thm. 7.2{]}}
\end{thm*}

\section{Fundamental Theorem\label{SEC - Fund Thm}}

\noindent The fundamental theorem of projective geometry {[}vSta47{]}
is required for many results, including Pascal's Theorem. Classically,
the fundamental theorem is derived from axioms of order and continuity.
For the  projective plane $\mathbb{P}$, since no axioms of order
and continuity have been adopted, the crucial component of the fundamental
theorem must be derived directly from an axiom: \\

\noindent \textbf{Axiom T.\label{Axiom T. }} \emph{If a projectivity
$\pi$ of a range or pencil onto itself has three distinct fixed elements,
then it is the identity $\iota$. }\\

Classically, Axiom T is often given the following equivalent form:
\emph{Let $\pi$ be a projectivity from a range onto itself, with
$\pi\neq\iota$, and distinct fixed points $M$ and $N$. If $Q$
is a point of the range distinct from both $M$ and $N$, then $Q^{\pi}\neq Q$.}
Constructively, this appears to be a stronger statement, since the
implication ``$\neg(Q^{\pi}=Q)$ implies $Q^{\pi}\neq Q$'' is constructively
invalid; thus we have: \\

\emph{Problem.} Give a proof of the apparently-stronger, alternative
statement for Axiom T, or show that it is constructively stronger.
\\

To prove that the alternative statement is constructively stronger
would require a Brouwerian counterexample. To determine the specific
nonconstructivities in a classical theory, and thereby to indicate
feasible directions for constructive work, Brouwerian counterexamples
are used, in conjunction with nonconstructive  omniscience principles.
A \emph{Brouwerian counterexample} is a proof that a given statement
implies an omniscience principle. In turn, an \emph{omniscience principle}
would imply solutions or significant information for a large number
of well-known unsolved problems. This method was introduced by L.
E. J. Brouwer {[}Brou08{]} to demonstrate that use of the law of excluded
middle inhibits mathematics from attaining its full significance.
A statement is considered \emph{constructively invalid} if it implies
an omniscience principle. The omniscience principles can be expressed
in terms of real numbers; the following are most often utilized: \\

\textit{\emph{Limited principle of omniscience (LPO). }}\emph{For
any real number $\alpha$, either $\alpha=0$ or $\alpha\neq0$. }\\

\textit{\emph{Weak limited principle of omniscience (WLPO).}} \emph{For
any real number $\alpha$, either $\alpha=0$ or $\neg(\alpha=0)$.
}\\

\textit{\emph{Lesser limited principle of omniscience (LLPO).}} \emph{For
any real number $\alpha$, either $\alpha\leq0$ or $\alpha\geq0$.
}\\
\emph{ }

\textit{\emph{Markov's principle. }}\emph{For any real number $\alpha$,
if $\neg(\alpha=0)$, then $\alpha\neq0$. }\\

For work according to Bishop-type strict constructivism, as followed
here, these principles, consequences of the law of excluded middle,
are used only to demonstrate the nonconstructive nature of certain
classical statements, and are not accepted for developing a constructive
theory. Markov's Principle, however, is used for work in recursive
function theory. 
\begin{thm*}
\emph{\label{Fundamental Theorem}}\textsc{\small{}Fundamental Theorem.}
Given any three distinct points $P$,$Q$,$R$ in a range $\overline{l}$,
and any three distinct points $P',Q',R'$ in a range $\overline{m}$,
there exists a unique projectivity $\pi:\overline{l}\rightarrow\overline{m}$
such that the points $P$,$Q$,$R$ map into the points $P',Q',R'$,
in the order given. \emph{{[}M16, Thm. 6.1{]}}
\end{thm*}
\begin{proof}
The existence of the required projectivity is provided by the second
theorem in Section \ref{SEC - Projectivities} above. Uniqueness,
however, requires Axiom T.
\end{proof}
\noindent Classically, the fundamental theorem is derived from axioms
of order and continuity; thus we have: \\

\emph{Problem.} Introduce constructive axioms of order and continuity
for the projective plane $\mathbb{P}$; derive Axiom T and the fundamental
theorem. \\

It follows from the fundamental theorem that any projectivity between
distinct ranges, or between distinct pencils, that has a fixed element
is a perspectivity {[}M16, Cor. 6.2{]}. A projectivity $\pi$ such
that $x^{\pi}\neq x$, for all elements $x$, is called \emph{nonperspective. }

The concept of projectivity is extended to the entire plane. A \emph{collineation
}of the projective  plane $\mathbb{P}$ is a bijection of the family
$\mathscr{P}$ of points, onto itself, that preserves collinearity
and noncollinearity. A collineation $\sigma$ induces an analogous
bijection $\sigma'$ of the family $\mathscr{L}$ of lines. A collineation
is \emph{projective} if it induces a projectivity on each range and
each pencil of the plane.

The following theorem is a constructivization of one of the main results
in the classical theory. 
\begin{thm*}
\label{Proj Coll Thm} A projective collineation with four distinct
fixed points, each three of which are noncollinear, is the identity.\emph{
{[}M16, Prop. 6.7{]}}
\end{thm*}
\begin{proof}
Let the collineation $\sigma$ have the fixed points $P,Q,R,S$ as
specified; thus the three distinct lines $PQ,PR,PS$ are fixed. The
mapping $\sigma'$ induces a projectivity on the pencil $P^{*}$;
by the fundamental theorem\emph{,} this projectivity is the identity.
Thus every line through $P$ is fixed under $\sigma'$; similarly,
the same is true for the other three points. 

Now let $X$ be any point on the plane. By three successive applications
of cotransitivity for points, we may assume that $X$ is distinct
from each of the points $P,Q,R$. Since $PQ\neq PR$, using cotransitivity
for lines we may assume that $XP\neq PQ$. Since $Q\neq P=XP\cdot PQ$,
it follows from Axiom C7 that $Q\notin XP$, and thus $XP\neq XQ$.
Since $X=XP\cdot XQ$, and the lines $XP$ and $XQ$ are fixed under
$\sigma'$, it follows that $\sigma X=X$. 
\end{proof}
\emph{Problem.} The above theorem ensures the uniqueness of a collineation
that maps four distinct points, each three of which are noncollinear,
into four distinct specified points, each three of which are also
noncollinear. Establish the existence of such a collineation for the
 projective plane $\mathbb{P}$. \\

The classical theory of the axis of homology\emph{ }has also been
constructivized\emph{. }
\begin{defn*}
\label{Axis Homology} Let $\pi:\overline{l}\rightarrow\overline{m}$
be a nonperspective projectivity between distinct ranges on the projective
plane $\mathbb{P}$. Set $O=l\cdot m$, $V=O^{\pi}$, and $U=O^{\pi^{-1}}$;
then the line $h=UV$ is called the \emph{axis of homology for} $\pi$.
{[}M16, Defn. 6.4{]}
\end{defn*}
The following theorem is the main result concerning the axis of homology;
the proof requires the fundamental theorem. 

\includegraphics[viewport=100bp 30bp 1266bp 500bp,clip,scale=0.4]{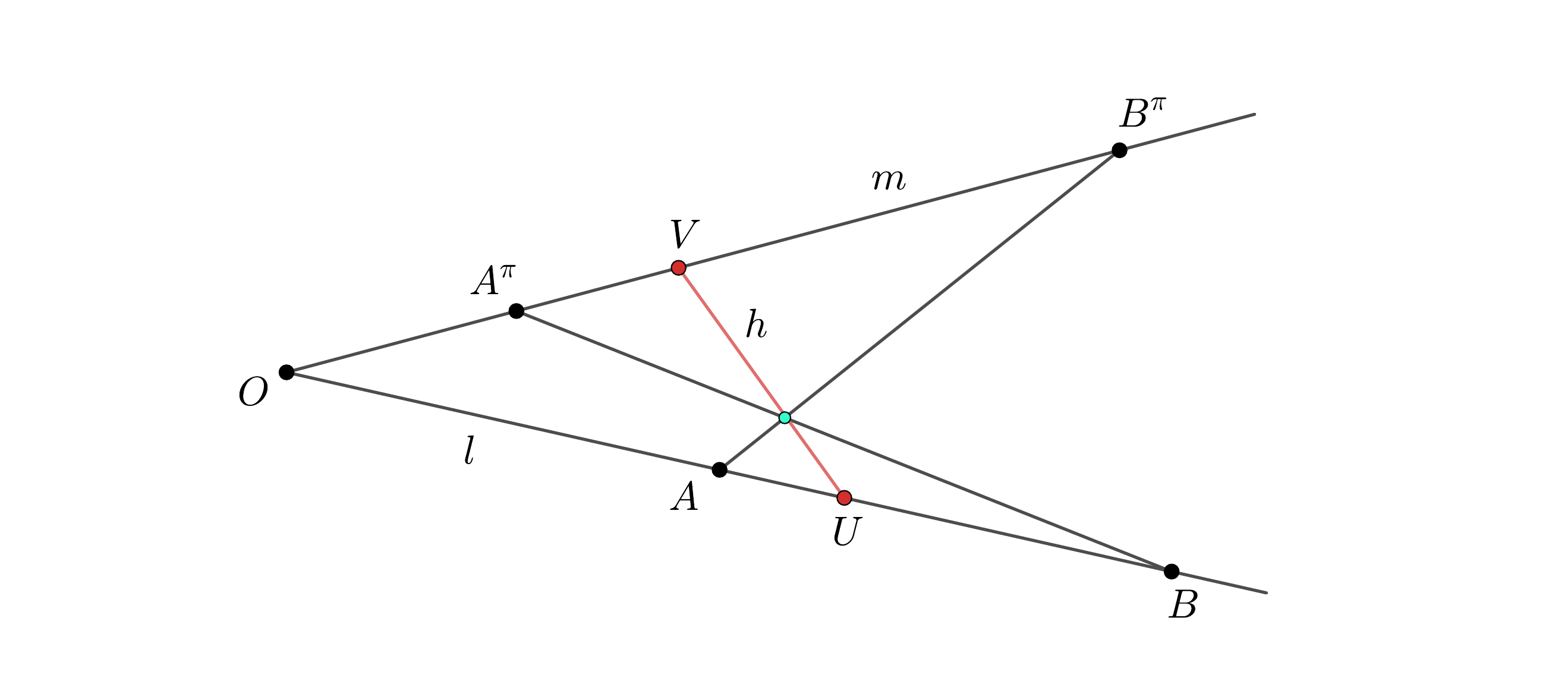}
\begin{thm*}
\label{Thm Axis Hom}Let $\pi:\overline{l}\rightarrow\overline{m}$
be a nonperspective projectivity between distinct ranges on the projective
plane $\mathbb{P}$. If $A$ and $B$ are distinct points on $l$,
each distinct from the common point $O$, then the point $AB^{\pi}\cdot BA^{\pi}$
lies on the axis of homology $h$. \emph{{[}M16, Thm. 6.5{]}}
\end{thm*}

\section{Conics\label{SEC - Conics}}

The conic sections have a long history; they were discovered by Menaechmus
(ca. 340 BC) and studied by the Greek geometers to the time of Pappus
of Alexandria (ca. 320 AD). The motivation for Menaechmus's discovery
was a geometrical problem, put forth by the oracle on the island of
Delos, the solution of which would have provided a remedy for the
Athenian plague of 430 BC. Unfortunately, Menaechmus's solution was
too late; see {[}Cox55, p. 79{]} for details. 

In the 17th century, an intense new interest in the conics arose in
connection with projective geometry. On a projective plane there is
no distinction between the hyperbola, parabola, and ellipse; these
arise only in the affine plane after a line at infinity is removed.
Which of the three forms results depends on whether that line meets
the conic at two, one, or no points. 

\subsection{Construction of a conic \label{Conic definition}}

\noindent Conics on the projective plane $\mathbb{P}$ are defined
by means of projectivities, using the method of Jakob Steiner (1796
\textendash{} 1863) {[}Ste32{]}. Alternatively, in classical works
conics are often defined by means of polarities, using the method
of von Staudt; see the problem stated at the end of Section \ref{SEC - Polarity}
below. 

\noindent \includegraphics[viewport=30bp 50bp 911bp 390bp,clip,scale=0.7]{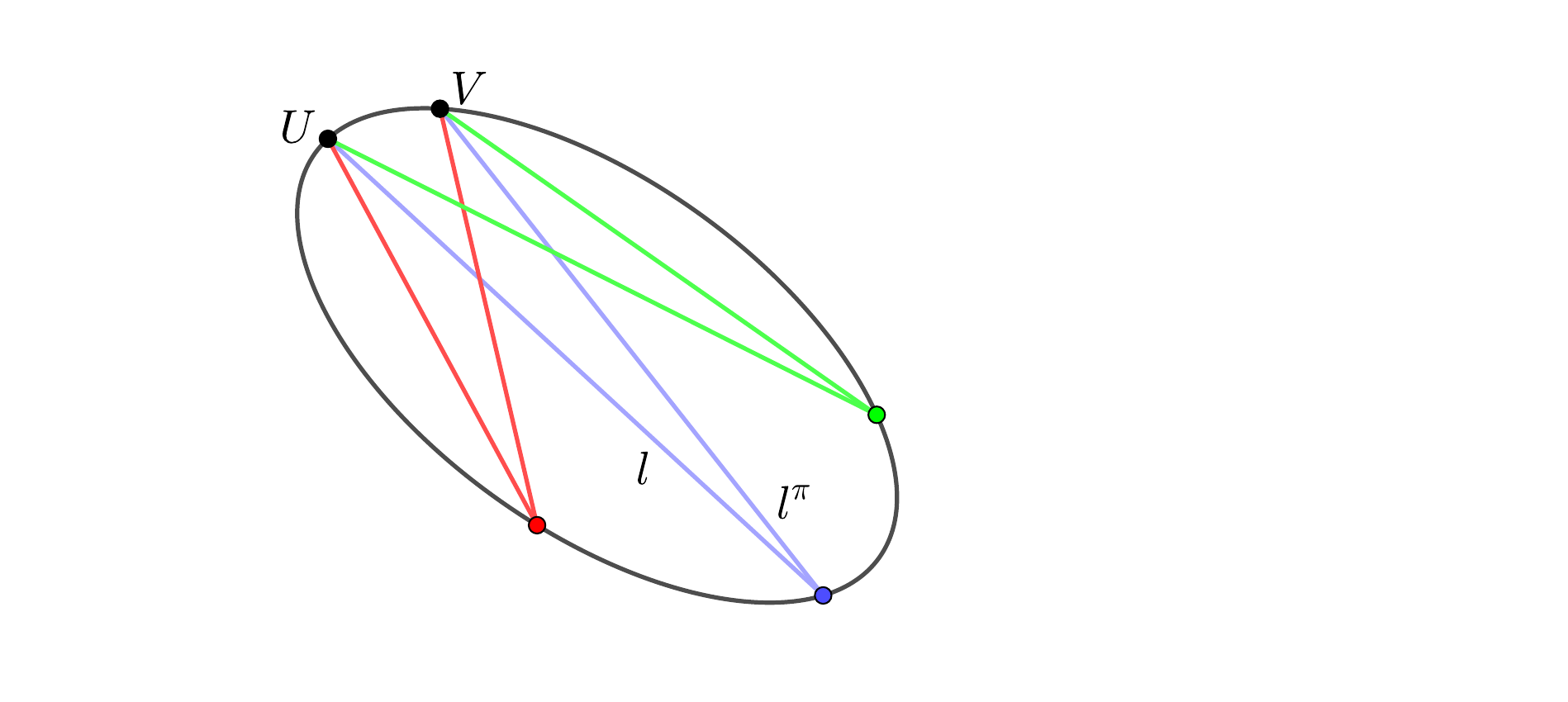}
\begin{defn*}
\label{Defn Conic}(Steiner) Let $\pi:U^{*}\rightarrow V^{*}$ be
a nonperspective projectivity between distinct pencils of lines on
the projective plane $\mathbb{P}$. The \emph{conic $\kappa=\kappa(\pi;U,V)$
defined by $\pi$ }is the locus of points $\{l\cdot l^{\pi}:l\in U^{*}\}$.
For any point $X$ on $\mathbb{P}$, we will say that \emph{$X$ lies
outside $\kappa$}, written $X\notin\kappa$, if $X\neq Y$ for all
points $Y$ on $\kappa$. {[}M16, Defn. 8.1{]} 

\medskip{}
\end{defn*}
\emph{Problem.} This definition, with the assumption that the given
projectivity is nonperspective, produces what is usually called a
\emph{non-singular }conic. Singular conics await constructive investigation.
\\

\subsection{Properties}

The next theorem establishes an essential property of a conic, an
analogue of the tightness property for inequalities; it can be viewed
as an extension of Axiom C6: ``If $\mbox{\ensuremath{\neg(P\notin l)}}$,
then $P\in l$.''
\begin{thm*}
Let \emph{$\kappa=\kappa(\pi;U,V)$}\textup{ }\textup{\emph{be a conic
}}on the  projective plane $\mathbb{P}$\textup{\emph{.}}\textup{
}For any point $X$ on $\mathbb{P}$\textup{,} if $\neg(X\notin\kappa)$,
then $X\in\kappa$. \emph{{[}M16, Prop. 8.2(d){]}}
\end{thm*}
\begin{proof}
Let $X$ be a point on the plane such that $\neg(X\notin\kappa)$.
By cotransitivity and symmetry, we may assume that $X\neq U$. Set
$z=UX$; then $Z=z\cdot z^{\pi}$ is a point of $\kappa$. Suppose
that $X\neq Z$. 

We now show that $X\neq Y$ for any point $Y$ of $\kappa$. Either
$Y\neq X$ or $Y\neq U$. We need to consider only the second case;
set $y=UY$, it follows that $Y=y\cdot y^{\pi}$. Either $Y\neq X$
or $Y\neq Z$; again, it suffices to consider the second case. Since
$Y\neq Z=z\cdot z^{\pi}$, it follows from Axiom C7 that either $Y\notin z$
or $Y\notin z^{\pi}$. In the first subcase, $y\neq z$. In the second
subcase, $y^{\pi}\neq z^{\pi}$, and since $\pi$ is a bijection we
again have $y\neq z$. Since $X\neq U=y\cdot z$, it follows that
$X\notin y$, and thus $X\neq Y$. 

The above shows that $X\notin\kappa$, contradicting the hypothesis.
It follows that $X=Z$, and hence $X\in\kappa.$ 
\end{proof}
Using this theorem and other preliminary results, many well-known
classical results are obtained constructively; for example, the following
basic result: 
\begin{thm*}
There exists a unique conic containing any given five distinct points,
each three of which are noncollinear. \emph{{[}M16, Prop. 8.3{]}}
\end{thm*}

\subsection{Pascal's Theorem}

Perhaps the most widely-known classical result concerning conics is
the following, due to Blaise Pascal (1623-1662)\emph{ }{[}Pas39{]}\emph{;}
it also has a constructive proof. 

\includegraphics[viewport=-30bp 0bp 688bp 400bp,clip,scale=0.6]{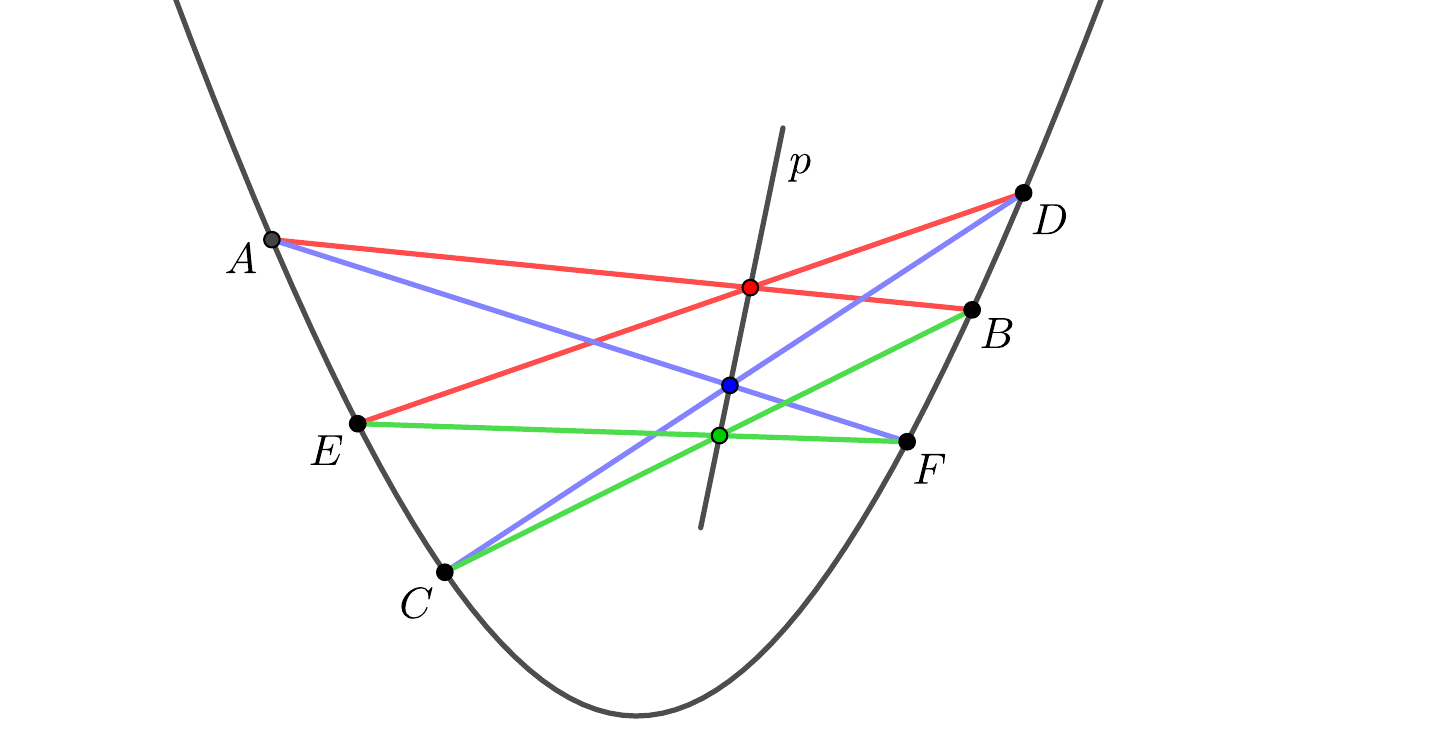}
\begin{thm*}
\textsc{\label{Pascal's-Theorem}}\textsc{\small{}Pascal's Theorem.}\emph{
}Let a simple hexagon $ABCDEF$ be inscribed in a conic $\kappa$.
Then the three points of intersection of the pairs of opposite sides
are distinct and collinear. \emph{{[}M16, Thm. 9.2{]} }
\end{thm*}
According to legend, Pascal gave in addition some four hundred corollaries.
Only one has been constructivized, it recalls a traditional construction
method for drawing a conic ``point by point'' on paper; for example,
as in {[}Y30, p. 68{]}. 

\includegraphics[viewport=-180bp 0bp 1222bp 530bp,clip,scale=0.52]{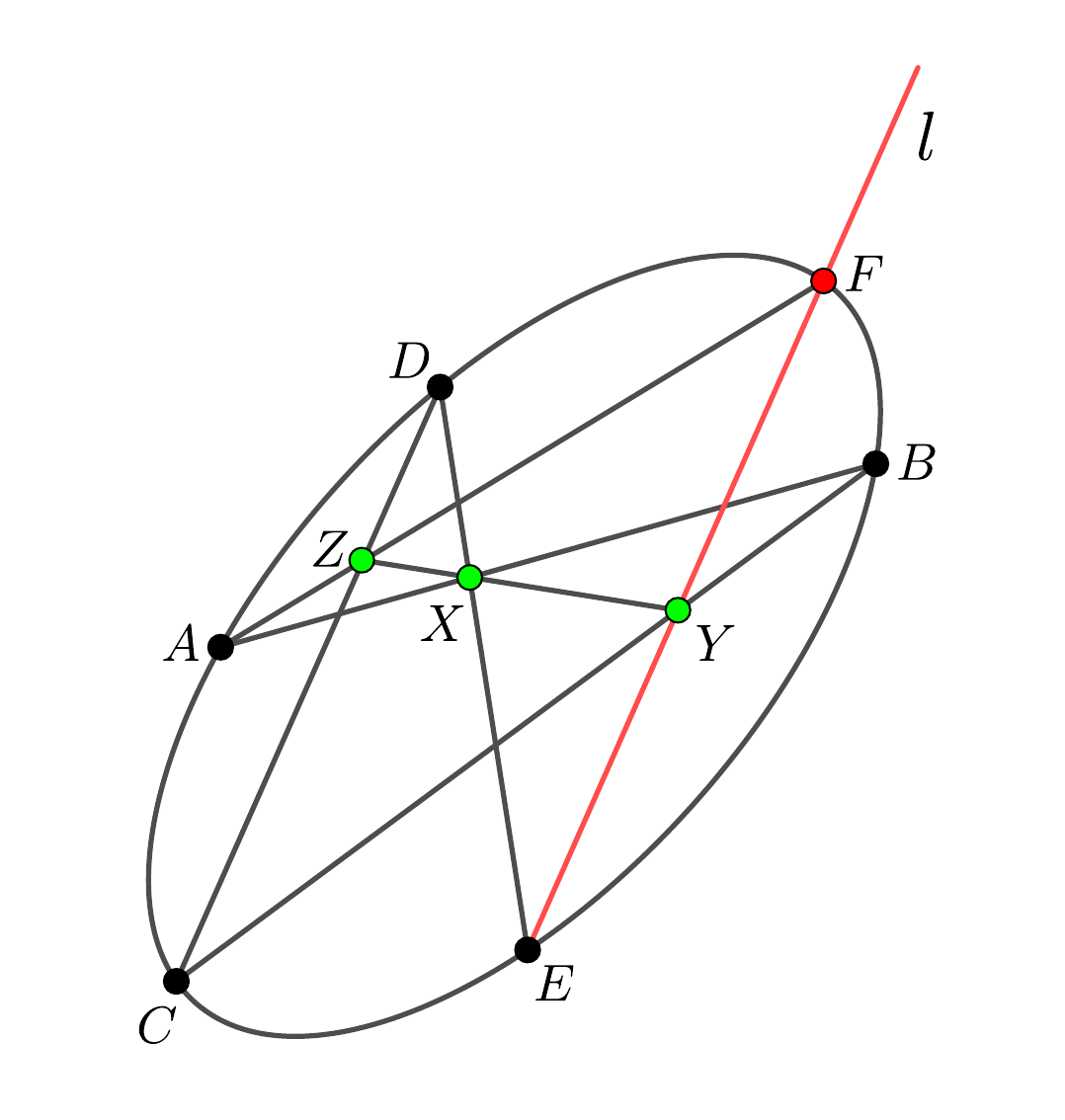}
\begin{cor*}
\label{Pascal Corollary}Let $A,B,C,D,E$ be five distinct points
of a conic $\kappa$. If $l$ is a line through $E$ that avoids each
of the other four points, and $l$ passes through a distinct sixth
point $F$ of $\kappa$, then
\[
F=l\cdot A(CD\cdot(AB\cdot DE)(BC\cdot l)).
\]

\noindent \emph{{[}M16, Cor. 9.3{]}}
\end{cor*}
\begin{proof}
The \emph{Pascal line} $p$ of the hexagon $ABCDEF$ passes through
the three distinct points $X=AB\cdot DE$, $Y=BC\cdot EF$, and $Z=CD\cdot AF$.
Since $A\notin CD$, we have $A\neq Z$, and it follows that $AF=AZ$.
Since $B\notin CD$, we have $BC\neq CD$, so by cotransitivity for
lines either $p\neq BC$ or $p\neq CD$. In the first case, since
$C\notin EF$, we have $C\neq Y=BC\cdot p$, and it follows from Axiom
C7 that $C\notin p$. Thus in both cases we have $CD\neq p$, and
$Z=CD\cdot p$. Finally, 
\begin{alignat*}{1}
F & =EF\cdot AF=l\cdot AZ=l\cdot A(CD\cdot p)\\
 & =l\cdot A(CD\cdot XY)=l\cdot A(CD\cdot(AB\cdot DE)(BC\cdot l))
\end{alignat*}
\end{proof}

\section{Polarity\label{SEC - Polarity}}

The role of symmetry in projective geometry reaches a peak of elegance
in the theory of polarity, introduced by von Staudt in 1847. 

A \emph{correlation} is a mapping of the points of the projective
plane to the lines, together with a mapping of the lines to the points,
that preserves collinearity and concurrence. The correlation is \emph{involutory}
if it is of order 2, and is then called a \emph{polarity.} A conic
determines a polarity; each point of the plane has a corresponding
polar, and each line has a corresponding pole. 

\subsection{Tangents and secants}

\noindent The construction of poles and polars determined by a conic
is dependent upon the existence of tangents and secants. A line $t$
that passes through a point $P$ on a conic \emph{$\kappa$} is said
to be \emph{tangent to $\kappa$ at $P$} if $P$ is the unique point
of $\kappa$ that lies on $t$. A line that passes through two distinct
points of a conic $\kappa$ is a \emph{secant} \emph{of} $\kappa$.
For the construction of poles and polars, it has been necessary to
adopt an additional axiom: \\

\noindent \textbf{Axiom P.\label{Axiom  P}} The tangents at any three
distinct points of a conic are nonconcurrent. \\

\emph{Problem. }Determine whether this axiom can be derived from the
others.  \\

The tangents and secants to a conic are related by means of projectivities;
the tangent at a point on a conic is the projective image of any secant
through the point:
\begin{thm*}
\label{tangent - secant - persp}Let $\kappa$ be a conic on the projective
plane $\mathbb{P}$, $P$ a point on $\kappa$, and $t$ a line passing
through $P$. The line $t$ is tangent to $\kappa$ at $P$ if and
only if for any point $Q$ of $\kappa$ with $Q\ne P$, if $s$ is
the secant $QP$, and $\pi$ is the nonperspective projectivity such
that $\kappa=\kappa(\pi;Q,P)$, then $t=s^{\pi}$. \emph{{[}M16, Prop.
10.2(b){]}}
\end{thm*}
This theorem ensures the existence of tangents. To establish the existence
of secants, it is first shown that a line through a point on a conic,
if not the tangent, is a secant:
\begin{lem*}
\label{Lemma - secant}Let $\kappa$ be a conic on the projective
plane $\mathbb{P}$, $P$ a point on $\kappa$, and $t$ the tangent
to $\kappa$ at $P$. If $l$ is a line passing through $P$, and
$l\neq t$, then $l$ passes through a second point $R$ of $\kappa$,
distinct from $P$; thus $l$ is a secant of $\kappa$. \emph{{[}M16,
Lm. 10.9{]}}
\end{lem*}
Using this lemma, the next theorem will provide the secants needed
for the following results. The need for this theorem contrasts with
complex geometry, where every line meets every conic.
\begin{thm*}
\label{Theorem - two secants}Let $\kappa$ be a conic on the projective
plane $\mathbb{P}$. Through any given point $P$ of the plane, at
least two distinct secants of $\kappa$ can be constructed. \emph{{[}M16,
Thm. 10.10(a){]}}
\end{thm*}
\begin{proof}
Select distinct points $A,B,C$ on $\kappa$, with tangents $a,b,c$.
By Axiom P, these tangents are nonconcurrent; thus the points $E=a\cdot b$
and $F=b\cdot c$ are distinct. Either $P\neq E$ or $P\neq F$; it
suffices to consider the first case. By Axiom C7, either $P\notin a$
or $P\notin b$. It suffices to consider the first subcase; thus $P\neq A$
and $PA\neq a$. It follows from the lemma that $PA$ is a secant. 

Denote the second point of $PA$ that lies on $\kappa$ by $R$, and
choose distinct points $A',B',C'$ on $\kappa$, each distinct from
both $A$ and $R$. With these three points, construct a secant through
$P$ using the above method; we may assume that it is $PA'$. Since
$A'\notin AR=PA$, it follows that $PA'\neq PA$. 
\end{proof}

\subsection{Construction of polars and poles}

The traditional method for defining a polar uses an inscribed quadrangle,
and must consider separately points on or outside a conic. Constructively,
this method is not available; thus polars are constructed by means
of harmonic conjugates. 

The discussion of harmonic conjugates in Section \ref{SEC - Harmonic-conjugates}
included an invariance theorem to show that the result of the construction
is independent of the selection of auxiliary elements. Now, the definition
of the polar of a point must be shown to be independent of the choice
of an auxiliary secant; the proof requires the invariance theorem
for harmonic conjugates. 

\includegraphics[viewport=70bp 0bp 1163bp 550bp,clip,scale=0.45]{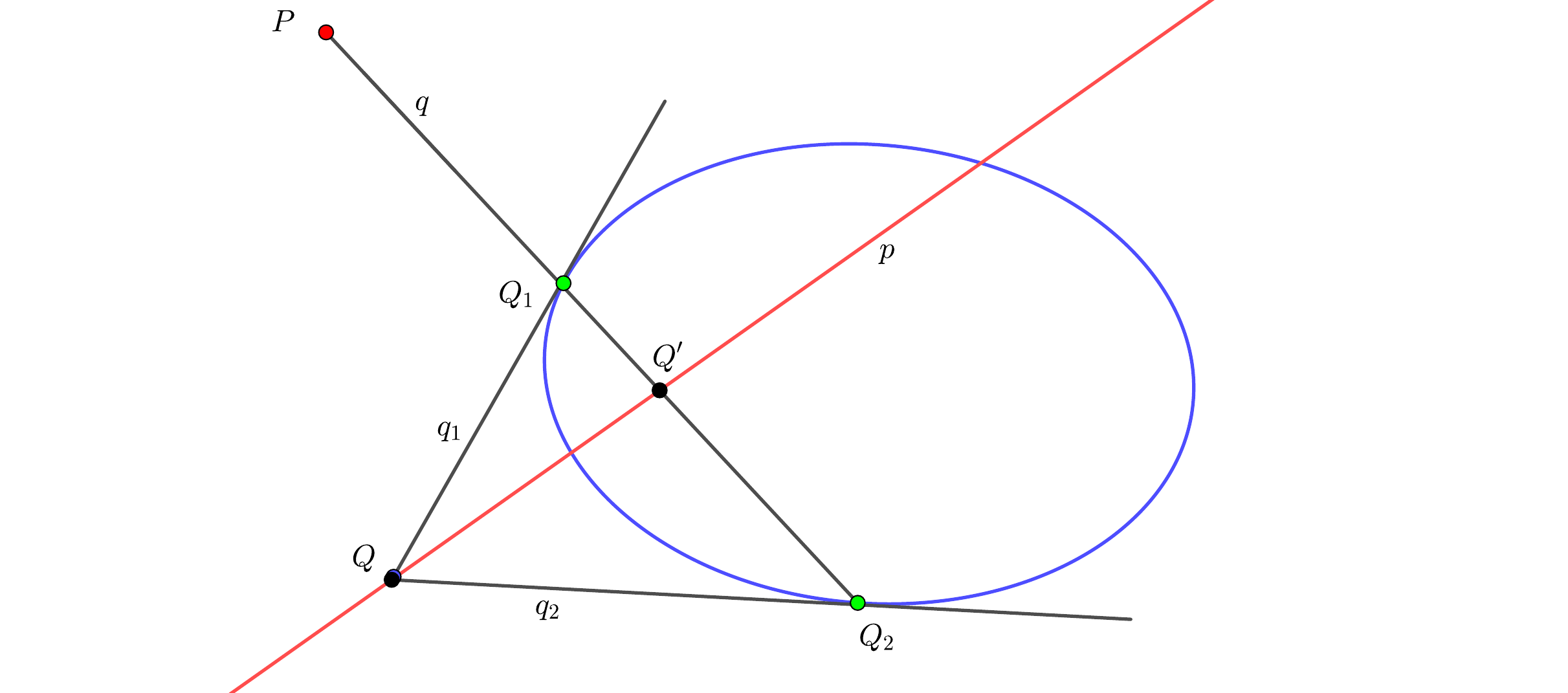}
\begin{thm*}
\textsc{\label{Polar construction}}\textsc{\small{}Construction of
a polar.} Let $\kappa$ be a conic on the projective plane $\mathbb{P}$,
and let $P$ be any point on the plane. Through the point $P$, construct
a secant $q$ of $\kappa$. Denote the intersections of $q$ with
$\kappa$ by $Q_{1}$ and $Q_{2}$, and let the tangents at these
points be denoted $q_{1}$ and $q_{2}$. Set $Q=q_{1}\cdot q_{2}$.
Set $Q'=h(Q_{1},Q_{2};P)$, the harmonic conjugate of $P$ with respect
to the base points $Q_{1},Q_{2}$. Then the line $p=QQ'$ is independent
of the choice of the secant $q$. \emph{{[}M16, Thm. 11.1{]}} 
\end{thm*}
\begin{defn*}
\label{Defn polar}Let $\kappa$ be a conic on the projective plane
$\mathbb{P}$, and let $P$ be any point on the plane. The line $p=QQ'$
in the above theorem is  called the \emph{polar of $P$ with respect
to $\kappa$. }{[}M16, Defn. 11.2{]}
\end{defn*}
Note that if $P$ lies on $\kappa$, then the polar of $P$ is the
tangent to $\kappa$ at $P$. The corollary below will relate this
constructive theory of polars to a classical construction that uses
quadrangles. The three \emph{diagonal points} of a quadrangle are
the intersection points of the three pairs of opposite sides. We adopt
Fano's Axiom: \emph{The diagonal points of any quadrangle are noncollinear.}
Gino Fano (1871-1952) studied finite projective planes, some of which
do not satisfy Fano's Axiom. 

\includegraphics[viewport=-150bp 0bp 1328bp 550bp,clip,scale=0.35]{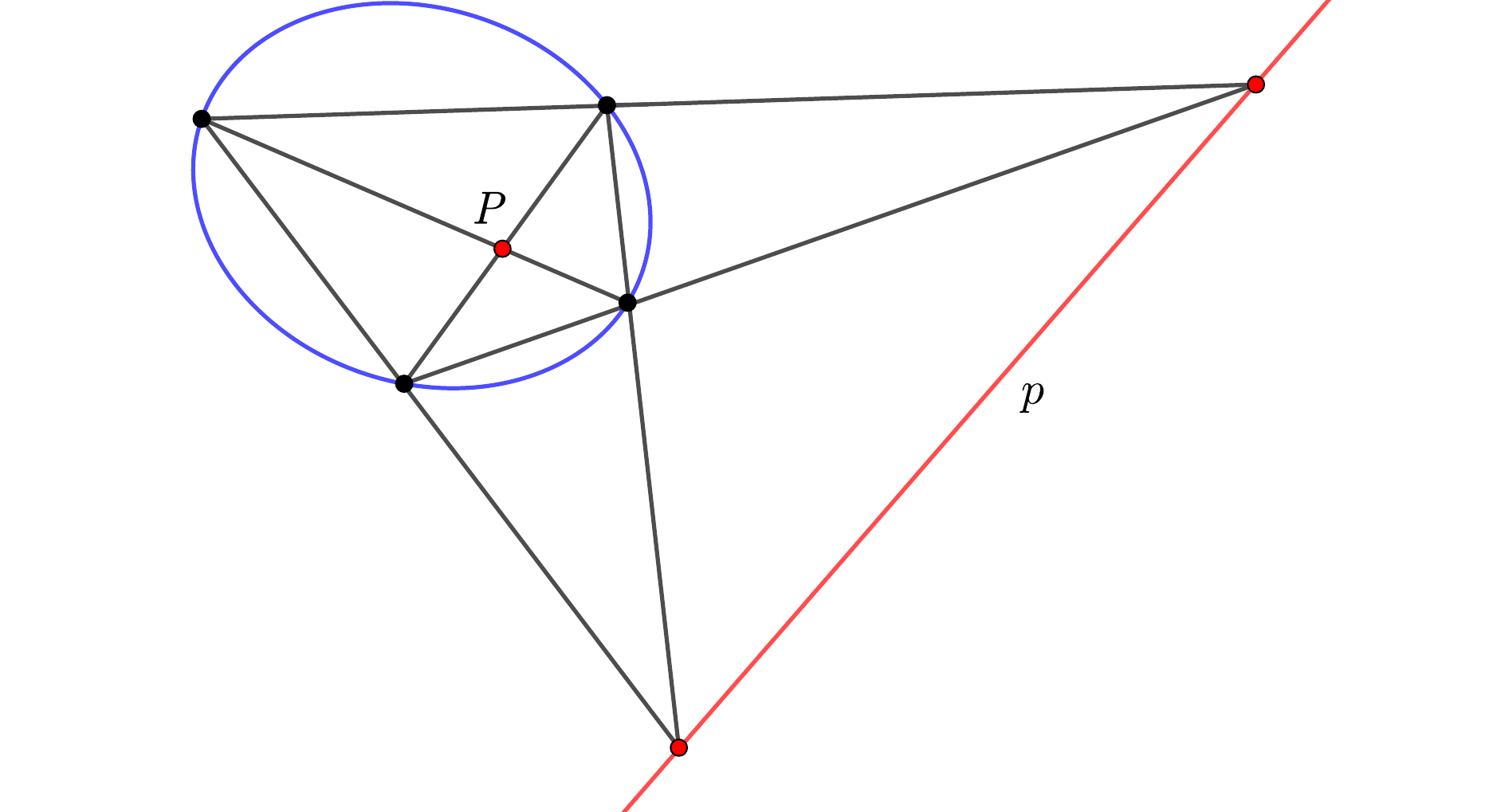}
\begin{cor*}
\label{Cor polar diagonal point}Let $\kappa$ be a conic on the projective
plane $\mathbb{P}$, and let $P$ be any point outside $\kappa$.
Inscribe a quadrangle in $\kappa$ with $P$ as one diagonal point.
Then the polar of $P$ is the line $p$ joining the other two diagonal
points. \emph{{[}M16, Cor. 11.4{]}}
\end{cor*}
Any three distinct points on a conic are noncollinear {[}M16, Prop.
8.2(b){]}. Thus, if $P$ is any diagonal point of a quadrangle inscribed
in a conic $\kappa$, it follows that $\neg(P\in\kappa)$. However,
it does not immediately follow that $P$ lies outside $\kappa$; thus
we have: \\

\emph{Problem. }If $\kappa$ is a conic on the projective plane $\mathbb{P}$,
and $P$ is a diagonal point of a quadrangle inscribed in $\kappa$,
show that  $P\notin\kappa$. \medskip{}

\begin{defn*}
\label{Defn pole}Let $\kappa$ be a conic on the projective plane
$\mathbb{P}$, and $l$ any line on $\mathbb{P}$. A construction
analogous to that of the above theorem results in a point $L$, called
the \emph{pole of} $l$ \emph{with respect to $\kappa$}. {[}M16,
Defn. 11.5{]}
\end{defn*}
The following theorem shows that any conic on the plane $\mathbb{P}$
determines a polarity:
\begin{thm*}
\label{polarity from conic} Let $\kappa$ be a conic on the projective
plane $\mathbb{P}$. If the line $p$ is the polar of the point $P$,
then $P$ is the pole of $p$, and conversely. \emph{{[}M16, Thm.
11.6(a){]}}
\end{thm*}
The definition of conic in Section \ref{SEC - Conics} used the Steiner
method {[}Ste32{]}, with projectivities. Later,  von Staudt {[}vSta47{]}
defined a conic by means of a polarity: a point lies on the conic
if its polar passes through the point. Classically, the two definitions
produce the same conics; thus we have: \\

\emph{Problem. }Construct correlations and polarities based on the
axioms for the projective plane $\mathbb{P}$, develop the theory
of conics constructively using the von Staudt definition, and prove
that von Staudt conics are equivalent to the Steiner conics constructed
in Section \ref{Conic definition} above. \\

\section{Consistency of the axiom system\label{SEC 8 Consistency}}

\noindent The consistency of the axiom system for the synthetic projective
plane $\mathbb{P}$ is established by an analytic model. A projective
plane $\mathbb{P}^{2}(\mathbb{R})$ is built from subspaces of the
linear space $\mathbb{R}^{3}$, using only constructive properties
of the real numbers. The axioms adopted for the synthetic plane $\mathbb{P}$
have been chosen to reflect the properties of the analytic plane $\mathbb{P}^{2}(\mathbb{R})$,
taking note of Bishop's thesis, ``All mathematics should have numerical
meaning'' {[}B67, p. ix; BB85, p. 3{]}. 

The model is built following well-known classical methods, adding
constructive refinements to the definitions and proofs. The analytic
plane $\mathbb{P}^{2}(\mathbb{R})$ consists of a family $\mathscr{P_{\mathrm{2}}}$
of points, and a family $\mathscr{L_{\mathrm{2}}}$ of lines; a \emph{point}
$P$ in $\mathscr{P_{\mathrm{2}}}$ is a subspace of dimension 1 of
the linear space $\mathbb{R}^{3}$, a \emph{line} $\lambda$ in $\mathscr{L_{\mathrm{2}}}$
is a subspace of dimension 2. The inequality relations, the incidence
relation, and the outside relation are defined by means of vector
operations. All the essential axioms adopted for the synthetic plane
$\mathbb{P}$, and all the required properties, such as cotransitivity,
tightness, and duality, are verified. 
\begin{thm*}
\label{Axioms verified}Axiom group C, and Axioms F, D, E, T, are
valid on the analytic projective plane $\mathbb{P}^{2}(\mathbb{R})$.
\emph{{[}M16, Thm. 14.2{]}}
\end{thm*}
The Brouwerian counterexample below shows that on the plane $\mathbb{P}^{2}(\mathbb{R})$
the validity of Axiom C3, which ensures the existence of a common
point for any two distinct lines, is dependent on the restriction
to distinct lines. By duality, the two statements of the example are
equivalent. The proof of the first statement is easier to visualize,
and can be described informally as follows: On $\mathbb{R}^{2}$,
thought of as a portion of $\mathbb{P}^{2}(\mathbb{R})$, consider
two points which are extremely near or at the origin, with $P$ on
the $x$-axis, and $Q$ on the $y$-axis. If $P$ is very slightly
off the origin, and $Q$ is at the origin, then the $x$-axis is the
required line. In the opposite situation, the $y$-axis would be required.
In any conceivable constructive routine, such a large change in the
output, resulting from a minuscule variation of the input, would reveal
a severe discontinuity, and a strong indication that the statement
in question is constructively invalid.
\begin{example*}
\noindent \label{Example lines common point}On the analytic projective
plane $\mathbb{P}^{2}(\mathbb{R})$, the following statements are
constructively invalid: 

(i)\emph{ Given any points $P$ and $Q$, there exists a line that
passes through both points.}

(ii)\emph{ Given any lines $\lambda$ and $\mu$, there exists a point
that lies on both lines. }

\noindent {[}M16, Example 14.1{]}
\end{example*}
\begin{proof}
By duality, it will suffice to consider the second statement. When
the non-zero vector $t=(t_{1},t_{2},t_{3})$ spans the point $T$
in $\mathbb{P}^{2}(\mathbb{R})$, we write $T=\langle t\rangle=\langle t_{1},t_{2},t_{3}\rangle$.
When the vectors $u,v$ span the line $\lambda$, and $w=u\times v$,
we write $\lambda=[w]=[w_{1},w_{2},w_{3}]$. The relation $T\in\lambda$
is defined by the inner product, $t\cdot w=0$. 

Let $\alpha$ be any real number, and set $\alpha^{+}=\max\{\alpha,0\}$,
and $\alpha^{-}=\max\{-\alpha,0\}$. Define lines $\lambda=[\alpha^{+},0,1]$
and $\mu=[0,\alpha^{-},1]$. By hypothesis, we have a point $T=\langle t\rangle=\langle t_{1},t_{2},t_{3}\rangle$
that lies on both lines. Thus $\alpha^{+}t_{1}+t_{3}=0$, and $\alpha^{-}t_{2}+t_{3}=0$.
If $t_{3}\neq0$, then we have both $\alpha^{+}\neq0$ and $\alpha^{-}\neq0$,
an absurdity; thus $t_{3}=0$. This leaves two cases. If $t_{1}\neq0$,
then $\alpha^{+}=0$, so $\alpha\leq0$, while if $t_{2}\neq0$, then
$\alpha^{-}=0$, so $\alpha\geq0$. Hence LLPO results.
\end{proof}
\emph{Problem. }Develop the theory of conics for the analytic plane
$\mathbb{P}^{2}(\mathbb{R})$; compare the results with those for
the synthetic plane $\mathbb{P}$. On the plane $\mathbb{P}^{2}(\mathbb{R})$,
determine the constructive validity of Axiom P of Section \ref{SEC - Polarity}
above. \\

\emph{Problem. }For the analytic projective plane $\mathbb{P}^{2}(\mathbb{R})$,
apply constructive methods to the study of harmonic conjugates, cross
ratios, and other topics of classical projective geometry. \\

\emph{Problem. }Although the model $\mathbb{P}^{2}(\mathbb{R})$ establishes
the consistency of the axiom system used for the  projective plane
$\mathbb{P}$, it remains to prove the independence of the axiom system,
or to reduce it to an independent system. 

\part{Projective extensions\label{PART II - Proj ext}}

\noindent The notion of infinity has mystified finite humans for millennia.
On the analytic projective plane $\mathbb{P}^{2}(\mathbb{R})$, where
points and lines are merely lines and planes through the origin in
$\mathbb{R}^{3}$, it is no surprise to notice that any two distinct
lines meet at a unique point. However, to envision two parallel lines
on $\mathbb{R}^{2}$ meeting at infinity requires some imagination.
Johannes Kepler (1571 \textendash{} 1630) invented the term ``focus''
in regard to ellipses, and stated that a parabola also has two foci,
with one at infinity. This idea was extended by Poncelet, leading
to the concepts of a line at infinity, and a projective plane. 

In the classical theory, a projective extension of an affine plane
is a fairly simple matter: each pencil of parallel lines determines
a point at infinity, at which the lines meet, and these points form
the line at infinity. A projective plane results, and the required
projective axioms are satisfied. The extension of the metric plane
$\mathbb{R}^{2}$ to a projective plane is often described heuristically,
with lamps and shadows; see, for example, {[}Cox55, Section 1.3{]}.

There have been at least three constructive attempts to extend an
affine plane to a projective plane. An extension by A. Heyting {[}H59{]}
uses elements called  ``projective points'' and ``projective lines''.
The extension constructed in {[}M14{]} uses elements called ``prime
pencils'' and ``virtual lines'', resulting in a projective plane
with different properties. The analytic projective plane $\mathbb{P}^{2}(\mathbb{R})$
discussed in Section \ref{SEC 8 Consistency} above, constructed using
subspaces of $\mathbb{R}^{3}$, can be viewed as an extension of the
metric plane $\mathbb{R}^{2}$; it also has distinctive properties. 

The differences between these several extensions involve the crucial
axiom concerning the existence of a point common to two lines, and
the cotransitivity property. The statement that any two distinct lines
have a common point is called the \emph{Common Point Property (CPP),}
while the \emph{Strong Common Point Property (SCPP)} is the same statement
without the restriction to distinct lines. The analytic extension
$\mathbb{P}^{2}(\mathbb{R})$ of $\mathbb{R}^{2}$ satisfies both
CPP\emph{ }and cotransitivity, but \emph{not} SCPP\emph{.} Neither
synthetic extension satisfies both cotransitivity and CPP\emph{.}
The Heyting extension satisfies cotransitivity, but the essential
axiom CPP\emph{ }has not been verified. On the virtual line extension,
CPP is satisfied, and even\emph{ }SCPP; however, cotransitivity is
constructively invalid, and this is now seen as a serious limitation.\emph{
}The analytic extension $\mathbb{P}^{2}(\mathbb{R})$ could be taken
as a standard; one might demand that the basic properties of $\mathbb{P}^{2}(\mathbb{R})$
hold in any acceptable synthetic extension, and then neither of the
two synthetic extensions would suffice. \\

\emph{Problem.} Construct a synthetic projective extension of an affine
plane which has the usual properties of a projective plane, including
both the common point property and cotransitivity. 

\section{Heyting extension\label{SEC - Heyting-.-.}}

In {[}H59{]}, A. Heyting adopts axioms for both affine and projective
geometry. Then, from a plane affine geometry $(\mathscr{P,L})$, Heyting
constructs an extension $(\Pi,\Lambda)$, consisting of \emph{projective
points} of the form 
\[
\mathfrak{P}(l,m)=\{n\in\mathscr{L}\,:\,n\cap l=l\cap m\,\,\,\text{or}\,\,\,n\cap m=l\cap m\}
\]

\noindent where $l,m\in\mathscr{L}$ with $l\neq m,$ and \emph{projective
lines }of the form
\[
\lambda(\mathfrak{A,B})=\{\mathfrak{Q}\in\Pi\,:\,\mathfrak{Q\cap A=A\cap B}\,\,\,\text{or\,\,\,}\mathfrak{Q\cap B=A\cap B}\}
\]

\noindent where $\mathfrak{A,B}\in\Pi$ with $\mathfrak{A\neq B.}$ 

With the Heyting definition of projective point, if the original two
lines $l$ and $m$ intersect, then $\mathfrak{P}(l,m)$ is the pencil
of all lines passing through the point of intersection, while if the
lines are parallel, then $\mathfrak{P}(l,m)$ is the pencil of all
lines parallel to the original two. In these cases, the definition
determines either a finite point of the extension, or a point on the
line at infinity. More significant is the fact that even when the
status of the two original lines is not known constructively, still
a projective point is (potentially) determined. Heyting comments on
the need for this provision as follows:
\begin{quotation}
. . . serious difficulties . . . are caused by the fact that not only
points at infinity must be adjoined to the affine plane, but also
points for which it is unknown whether they are at infinity or not.
{[}H59, p. 161{]}
\end{quotation}
A projective line is determined by two distinct projective points.
The definition is based on the lines common to the two projective
points; i.e., the lines common to two pencils of lines. For example,
in the simplest case, if the two projective points are finite, then
these are the pencils of lines through distinct points in the original
affine plane, and there is a single common line, connecting these
finite points, of which the projective line is an extension. In the
case of two distinct pencils of parallel lines; the pencils have no
common line, each determines a point at infinity, and the resulting
projective line is the line at infinity. Again, even when the status
of the original projective points is not known constructively, still
a projective line is determined. The distinctive, and perhaps limiting,
features of the Heyting extension are the requirements that the construction
of a projective point depends on a given pair of distinct finite lines,
and the construction of a projective line depends on a pair of distinct
projective points previously constructed. 

Nearly all the axioms for a projective plane are then verified, although
the most essential axiom, which states that two distinct lines have
a common point, escapes proof. The axiom considered in {[}H59{]} is
the weaker version, designated above as the common point property,
CPP, involving distinct lines. In {[}M13{]}, Heyting's axioms for
affine geometry are verified for $\mathbb{R}^{2}$, and a Brouwerian
counterexample is given for the Heyting extension, showing that the
stronger form of the axiom, SCPP\emph{, }involving arbitrary lines,
is constructively invalid, with the following attempted justification: 
\begin{quotation}
This counterexample concerns the full common point axiom, rather than
the limited Axiom P3 as stated in {[}H59{]}, where only distinct lines
are considered. An investigation into the full axiom is necessary
for a constructive study based upon numerical meaning, as proposed
by Bishop. Questions of distinctness are at the core of constructive
problems; any attempted projective extension of the real plane is
certain to contain innumerable pairs of lines which may or may not
be distinct. {[}M13, p. 113{]} 
\end{quotation}
However, taking note of the analytic model $\mathbb{P}^{2}(\mathbb{R})$,
for which CPP is verified, but SCPP is constructively invalid, CPP
now appears as a reasonable goal for an extension; thus we have: \\

\emph{Problem. }Complete the study of the projective extension of
{[}H59{]}; verify Heyting's Axiom P3 (CPP), or construct a Brouwerian
counterexample\emph{. }

\section{Virtual line extension\label{SEC Virtual line extension} }

Any attempt to build a constructive projective extension of an affine
plane encounters difficulties due to the indeterminate nature of arbitrary
pencils of lines. Classically, a pencil of lines is either the family
of lines passing through a given point, or a family of parallel lines.
An example of a family of lines is easily formed from two lines which
might be distinct, intersecting or parallel, or might be identical.
To obtain the strong common point property, SCPP, in a constructive
projective extension, the corresponding pencil must include both these
lines, so that it will determine a point of the extension common to
both extended lines, whether distinct or not. Thus the definition
of pencil must not depend upon a pair of lines previously known to
be distinct. 

In the projective extension of {[}M14{]}, the definition of pencil
is further generalized; rather than depending upon specific finite
lines, it involves the intrinsic properties of a family of lines.
Included are pencils of unknown type, with non-specific properties,
and pencils for which no lines are known to have been previously constructed. 

The definition of line in the extension is independent of the definition
of point; it will depend directly upon a class of generalized lines
in the finite plane, called virtual lines. 

\subsection{Definition\label{VL Ext Defn}}

The virtual line extension of {[}M14{]} is based on a incidence plane
$\mathscr{G=}(\mathscr{P,L}),$ consisting of a family $\mathscr{P}$
of points and a family $\mathscr{L}$ of lines, with constructive
axioms, definitions, conventions, and results as delineated in {[}M07{]}. 
\begin{defn*}
\label{Defn Q* l* Pencil }Let $\mathscr{G=}(\mathscr{P,L})$ be an
incidence plane.

\noindent $\bullet$ For any point $Q\in\mathscr{P},$ define 

\noindent 
\[
Q^{*}=\{l\in\mathscr{L}:Q\in l\}.
\]

\noindent $\bullet$ For any line $l\in\mathscr{L},$ define 

\noindent 
\[
l^{*}=\{m\in\mathscr{L}:m\parallel l\}.
\]

\noindent $\bullet$ A family of lines  of the form $Q^{*}$, or of
the form $l^{*}$, is called a \emph{regular} \emph{pencil}. 

\noindent $\bullet$ A family of lines $\alpha$ is called a \emph{pencil}
if it contains no fewer than two lines, and satisfies the following
condition: \emph{If $l$ and $m$ are distinct lines in $\alpha$
with $l,m\in\rho,$ where $\rho$ is a regular pencil, then $\alpha\subset\rho$. }

\noindent $\bullet$ A pencil of the form $Q^{*}$ is called a \emph{point
pencil}. 

\noindent $\bullet$ A pencil $\alpha$ with the property that $l\parallel m,$
for any lines $l$ and $m$ in $\alpha,$ is called a \emph{parallel
pencil}. 

\noindent {[}M14, Defn. 2.1{]}
\end{defn*}
In the extension, a point pencil $Q^{*}$, consisting of all lines
through $Q$, will represent the original finite point $Q$. A pencil
$l^{*}$, consisting of all lines parallel to the line $l$, will
result in an infinite point. However, the extension also admits parallel
pencils which need not arise from given lines, but which nevertheless
result in points at infinity. 

\subsection{Virtual lines\label{Virtual lines}}

A problem that arises in the construction of a projective extension
is the difficulty in determining the nature of an arbitrary line in
the extension, by means of an object in the original plane. If a line
$\lambda$  on the extended plane contains a finite point, then the
set $\lambda_{f}$, of all finite points on $\lambda$, is a line
in the original plane. However, if $\lambda$ is the line at infinity,
then $\lambda_{f}$ is void. Since constructively it is in general
not known which is the case, we adopt the following: 
\begin{defn*}
\noindent \label{Defn Virtual line}A set $p$ of points in $\mathscr{P}$
is said to be a \emph{virtual line} if it satisfies the following
condition: \emph{If $p$ is inhabited, then $p$ is a line. }{[}M14,
Defn. 3.1{]} 
\end{defn*}
Given any virtual lines $p$ and $q$, one can  construct a pencil
$\varphi(p,q)$ that contains each of the virtual lines $p$ and $q$,
if it is a line {[}M14, Thm. 3.4{]}.

The notion of virtual line also helps in resolving a problem that
arises in connection with pencils of lines. The family of lines common
to two distinct pencils may consist of a single line (as in the case
of two point pencils, or a point pencil and a regular parallel pencil),
or it may be void (as in the case of two regular parallel pencils);
constructively, it is in general unknown which alternative holds.
The following definition provides a tool for dealing with this situation. 
\begin{defn*}
\label{Defn Core}For any distinct pencils $\alpha$ and $\beta,$
define 

\[
\alpha\sqcap\beta=\{Q\in\mathscr{P}:Q\in l\in\alpha\cap\beta\text{ for some line }l\in\mathscr{L}\}.
\]

\noindent The set of points $\alpha\sqcap\beta$ is called the \emph{core
}of the pair of pencils $\alpha,\beta$. {[}M14, Defn. 3.2{]}
\end{defn*}
The core, as a set of finite points (which might  be void), is a constructive
substitute for a possible line that is common to two pencils. 
\begin{thm*}
\label{Thm Core is virtual line}For any distinct pencils $\alpha$
and $\beta$, the core \textup{$\alpha\sqcap\beta$} \textup{is a
}virtual\textup{\emph{ line}}\textup{. }\emph{{[}M14, Lm. 3.3{]}}
\end{thm*}

\subsection{Extension\label{Extension - virtual}}

Points of the extension, called \emph{e-points}, are defined using
a selected class of pencils of lines, called \emph{prime pencils;}
the prime pencil $\alpha$ determines the e-point $\overline{\alpha}$.
Lines in the extension are not formed from previously constructed
e-points; they are direct extensions of objects in the original plane.
Lines of the extension, called \emph{e-lines}, are defined using a
selected class of virtual lines, called \emph{prime virtual lines;}
the prime virtual line $p$ in the finite plane determines the e-line
$\lambda_{p}$ in the extended plane. 

The projective plane $\mathscr{G^{*}}=(\mathscr{P^{*}},\mathscr{L^{*}})$,
where $\mathscr{P^{*}}$ is the family of e-points, and $\mathscr{L^{*}}$
is the family of e-lines, is the \emph{projective extension} of the
incidence plane $\mathscr{G=}(\mathscr{P,L})$. The axioms of projective
geometry are verified for the extension. The following theorems are
the main results; the proof outlines will exhibit the symmetry of
the construction, and the utility of adopting independent definitions
for e-points and e-lines. 
\begin{thm*}
\label{Thm. Two e-points}On the projective extension $\mathscr{G^{*}}$
of the plane $\mathscr{G}$, there exists a unique e-line passing
through any two distinct e-points. \emph{{[}M14, Thm. 5.3{]}}
\end{thm*}
\noindent \emph{Proof outline. }The given e-points $\overline{\alpha}$
and $\overline{\beta}$ originate from pencils $\alpha$ and $\beta$;
the core $p=\alpha\sqcap\beta$ of these pencils is a virtual line
on the finite plane. This virtual line $p$ determines an e-line $\lambda_{p}$
in the extension, which passes through both e-points $\overline{\alpha}$
and $\overline{\beta}$. 
\begin{thm*}
\label{Thm Two e-lines}On the projective extension $\mathscr{G^{*}}$
of the plane $\mathscr{G}$, any two e-lines have an e-point in common.
If the e-lines are distinct, then the common e-point is unique. \emph{{[}M14,
Thm. 5.5{]}}
\end{thm*}
\noindent \emph{Proof outline. }The given e-lines $\lambda_{p}$ and
$\lambda_{q}$ originate from virtual lines $p$ and $q$; these virtual
lines determine a pencil $\gamma=\varphi(p,q)$ of lines on the finite
plane. This pencil $\gamma$ determines an e-point $\overline{\gamma}$
in the extension, which lies on both e-lines $\lambda_{p}$ and $\lambda_{q}$.
\\

Several definitions in {[}M14{]} involve negativistic concepts; for
example, Definition 2.1 for pencil, and Definition 3.1 for distinct
virtual lines.  Can this be avoided? Generally in constructive mathematics
one tries to avoid negativistic concepts, but perhaps some are unavoidable
in constructive geometry; thus we have: \\

\emph{Problem. }Modify the virtual line extension so as to avoid negativistic
concepts as far as possible. 

\subsection{The cotransitivity problem\label{Subsec: cotransitivity}}

There is what might be called an irregularity of the extension plane
$\mathscr{G^{*}}$, the constructive invalidity of cotransitivity;
this is revealed by a Brouwerian counterexample: 
\begin{example*}
\label{BrCtEx Cotransi cxtvinv}On the virtual line projective extension
of the plane $\mathbb{R}^{2}$, the cotransitivity property for e-points
is constructively invalid. {[}M14, p. 705{]}
\end{example*}
\begin{proof}
\noindent Given any real number $c$, construct the virtual line 
\[
p=\{(t,0):t\in\mathbb{R}\,\,\text{and}\,\,c=0\}\cup\{(0,t):t\in\mathbb{R}\,\,\text{and}\,\,c\ne0\}
\]
and consider the e-point $\overline{\gamma}$ determined by the pencil
$\gamma=\varphi(p,p)$. 

Let the $x$-axis be denoted by $l_{0}$; the pencil $l_{0}^{*}$
of horizontal lines then determines the e-point $\overline{l_{0}^{*}}$.
Similarly, we have the $y$-axis $m_{0}$, the pencil $m_{0}^{*}$
of vertical lines, and the e-point $\overline{m_{0}^{*}}$. By hypothesis,
$\overline{\gamma}\neq\overline{l_{0}^{*}}$ or $\overline{\gamma}\neq\overline{m_{0}^{*}}$.
In the first case, suppose that $c=0$. Then $p$ is the $x$-axis
and $\overline{\gamma}=\overline{l_{0}^{*}},$ a contradiction; thus
we have $\neg(c=0)$. In the second case, we find that $c=0$. Hence
WLPO results. 
\end{proof}
\medskip{}

\emph{Problem.} Modify the virtual line extension, so that the common
point property and cotransitivity are both valid. It is then likely
that the strong common point property will not be valid; in that case,
provide a Brouwerian counterexample.

\section{Analytic extension\label{SEC:-analytic extension} }

The analytic projective plane $\mathbb{P}^{2}(\mathbb{R})$ described
in Section \ref{SEC 8 Consistency} is constructed from subspaces
of the linear space $\mathbb{R}^{3}$, using only constructive properties
of the real numbers. This projective plane can be viewed as an extension
of the affine plane $\mathbb{R}^{2}$. 

\includegraphics[viewport=40bp -80bp 915bp 350bp,clip,scale=0.65]{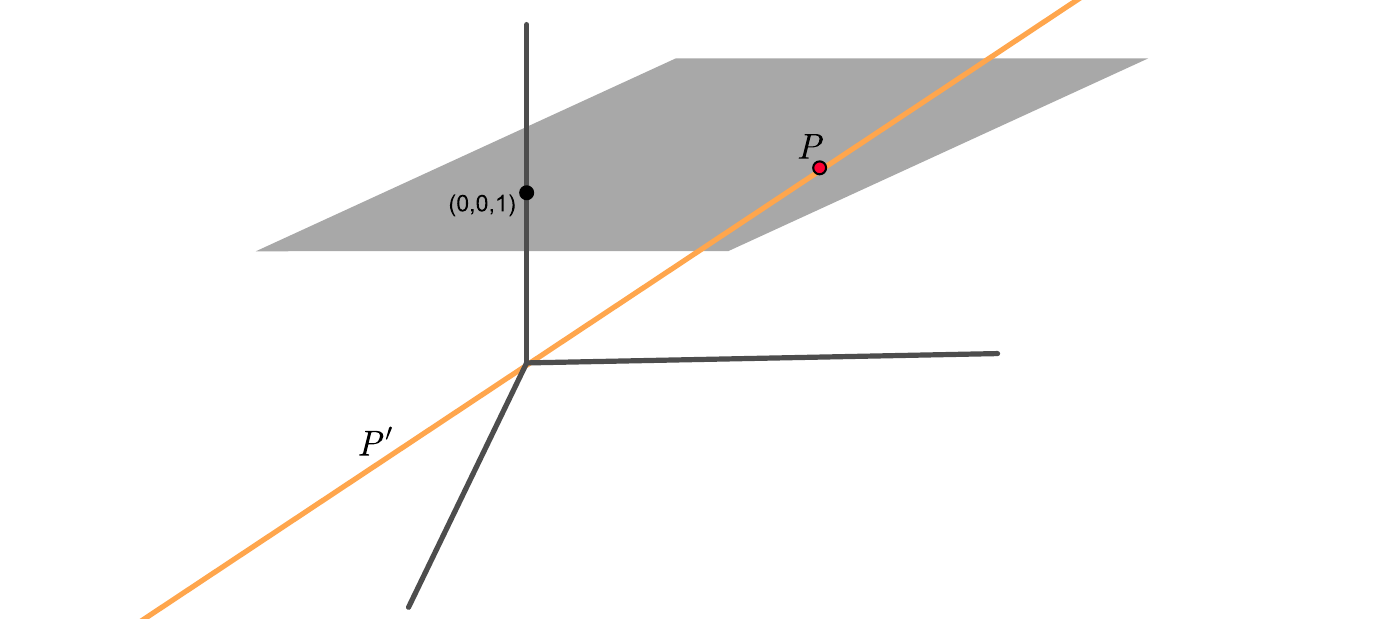}

The plane $z=1$ in $\mathbb{R}^{3}$ is viewed as a copy of $\mathbb{R}^{2}$.
A point $P$ on the plane $z=1$ corresponds to the point $P'$ of
the extension $\mathbb{P}^{2}(\mathbb{R})$ that, as a line through
the origin in $\mathbb{R}^{3}$, contains $P$. A point of $\mathbb{P}^{2}(\mathbb{R})$,
that is a horizontal line through the origin in $\mathbb{R}^{3}$,
is an infinite point of the extension. A line $l$ on the plane $z=1$
corresponds to the line $l'$ of $\mathbb{P}^{2}(\mathbb{R})$ that,
as a plane through the origin in $\mathbb{R}^{3}$, contains $l$.
The line of intersection of this plane with the $xy$-plane is the
point at infinity on $l'$. In this way, $\mathbb{P}^{2}(\mathbb{R})$
is seen as a projective extension of $\mathbb{R}^{2}$, with the \emph{$xy$-}plane
as the line at infinity.

The plane $\mathbb{P}^{2}(\mathbb{R})$ satisfies both the common
point property and cotransitivity. However, as a projective extension
of the specific plane $\mathbb{R}^{2}$, it does not provide an extension
of an arbitrary affine plane; thus we have: \\

\emph{Problem. }Construct a synthetic projective extension of an arbitrary
affine plane, having both the common point property and the cotransitivity
property. 

\section*{{\normalsize{}References \label{References}}}

~~~~{[}B65{]} E. Bishop, Book Review: \emph{The Foundations of
Intuitionistic Mathematics, }by S. C. Kleene and R. E. Vesley, \emph{Bull.
Amer. Math. Soc.} 71:850-852, 1965. 

{[}B67{]} E. Bishop, \emph{Foundations of Constructive Analysis,}
McGraw-Hill, New York, 1967. 

{[}B73{]} E. Bishop, \emph{Schizophrenia in Contemporary Mathematics},
AMS Colloquium Lectures, Missoula, Montana, Amer. Math. Soc., Providence,
RI, 1973; reprint in \emph{Contemporary Mathematics} 39:l-32, 1985. 

{[}B75{]} E. Bishop, The crisis in contemporary mathematics, \emph{Proceedings
of the American Academy Workshop on the Evolution of Modern Mathematics},
Boston, 1974\emph{; Historia Math.} 2:507\textendash 517, 1975.

{[}BB85{]} E. Bishop and D. Bridges, \emph{Constructive Analysis,}
Springer-Verlag, Berlin, 1985.

{[}Be10{]} M. Beeson, Constructive geometry, \emph{10th Asian Logic
Conference,} World Sci. Publ., Hackensack, NJ, 2010, pp. 19\textendash 84.

{[}Be16{]} M. Beeson, Constructive geometry and the parallel postulate,
\emph{Bull. Symb. Log.} 22:1\textendash 104, 2016. 

{[}BR87{]} D. Bridges and F. Richman, \emph{Varieties of Constructive
Mathematics,} Cambridge University Press, Cambridge, UK, 1987. 

{[}Brou08{]} L. E. J. Brouwer, De onbetrouwbaarheid der logische principes,
\emph{Tijdschrift voor Wijsbegeerte} 2:152-158, 1908; translation,
The unreliability of the logical principles, A. Heyting (ed.), \emph{L.
E. J. Brouwer: Collected Works 1: Philosophy and Foundations of Mathematics,}
Elsevier, Amsterdam-New York, 1975, pp. 107\textendash 111. 

{[}Brou24{]} L. E. J. Brouwer, Intuitionistische Zerlegung mathematischer
Grundbegriffe; reprint, \emph{L. Brouwer, Collected Works, Vol. 1,}
North-Holland, Amsterdam, 1975, pp. 275\textendash 280.

{[}BV06{]} D. Bridges and L. Vî\c{t}\u{a},\textsc{ }\emph{Techniques
of Constructive Analysis,} Springer, New York, 2006.

{[}Cox55{]} H. S. M. Coxeter, \emph{The Real Projective Plane, 2nd
ed.,} Cambridge University Press, Cambridge, UK, 1955; reprint, Lowe
and Brydone, London, 1960. 

{[}Cre7\textsc{3{]} }L. Cremona,\textsc{ }\emph{Elementi di geometria
projettiva, }G. B. Paravia e Comp., Torino, 1873; translation, C.
Leudesdorf, \emph{Elements of Projective Geometry,} Clarendon Press,
Oxford, 1985; reprint, Forgotten Books, Hong Kong, 2012. 

{[}Des64{]} G. Desargues, \emph{Oeuvres De Desargues,} Leiber, 1864.

{[}H28{]} A. Heyting, Zur intuitionistischen Axiomatik der projektiven
Geometrie, \emph{Math. Ann. }98:491-538, 1928. 

{[}H59{]} A. Heyting, Axioms for intuitionistic plane affine geometry,
L. Henkin, P. Suppes, A. Tarski (eds.), \emph{The Axiomatic Method,
with special reference to geometry and physics: Proceedings of an
international symposium held at the University of California, Berkeley,
December 26, 1957 - January 4, 1958,} North-Holland, Amsterdam, 1959,
pp. 160-173. 

{[}H66{]} A. Heyting, \emph{Intuitionism: An Introduction,} North-Holland,
Amsterdam, 1966. 

{[}Kle72{]} F. Klein, \emph{Vergleichende Betrachtungen über neuere
geometrische Forschungen,} A. Deichert, Erlangen, 1872; translation,
M. W. Haskell, A comparative review of researches in geometry, \emph{Bull.
New York Math. Soc.} 2:215-249, 1893. 

{[}Leh17{]} D. N. Lehmer, \emph{An Elementary Course in Synthetic
Projective Geometry,} Ginn, Boston, 1917. 

{[}LV98{]} M. Lombard and R. Vesley, A common axiom set for classical
and intuitionistic plane geometry, \emph{Ann. Pure Appl. Logic} 95:229-255,
1998.

{[}M07{]} M. Mandelkern, Constructive coördinatization of Desarguesian
planes, \emph{Beiträge Algebra Geom.} 48:547-589, 2007.

{[}M13{]} M. Mandelkern, The common point problem in constructive
projective geometry, \emph{Indag. Math.} \emph{(N.S.)} 24:111-114,
2013. 

{[}M14{]} M. Mandelkern, Constructive projective extension of an incidence
plane, \emph{Trans. Amer. Math. Soc. }366:691-706, 2014. 

{[}M16{]} M. Mandelkern, A constructive real projective plane, \emph{J.
Geom. }107:19-60, 2016. 

{[}M18{]} M. Mandelkern, Constructive harmonic conjugates, \emph{Beiträge
Algebra Geom.} 60:391\textendash 398, 2019. 

{[}Pam98{]} V. Pambuccian, Zur konstruktiven Geometrie Euklidischer
Ebenen, \emph{Abh. Math. Sem. Univ. Hamburg} 68:7\textendash 16, 1998. 

{[}Pam01{]} V. Pambuccian, Constructive axiomatization of plane hyperbolic
geometry, \emph{MLQ Math. Log. Q.} 47:475\textendash 488, 2001. 

{[}Pam03{]} V. Pambuccian, Constructive axiomatization of non-elliptic
metric planes, \emph{Bull. Polish Acad. Sci. Math.} 51:49\textendash 57,
2003. 

{[}Pam05{]} V. Pambuccian, Elementary axiomatizations of projective
space and of its associated Grassmann space, \emph{Note Mat.} 24:129\textendash 141,
2005. 

{[}Pam11{]} V. Pambuccian, The axiomatics of ordered geometry; I.
Ordered incidence spaces, \emph{Expo. Math.} 29:24\textendash 66,
2011. 

{[}Pas39{]} B. Pascal, \emph{Essai pour les coniques,} 1639. 

{[}Pic75{]} G. Pickert, \emph{Projektive Ebenen},\emph{ 2. Aufl.,}
Springer-Verlag, Berlin-New York, 1975. 

{[}Pon22{]} J-V. Poncelet, \emph{Traité des Propriétés Projectives
des Figures,} Gauthier-Villars, Paris, 1822. 

{[}Ste32{]} J. Steiner, \emph{Sytematische Entwickelung der Abhängigkeit
geometrischer Gestalten von einander,} G. Fincke, Berlin, 1832. 

{[}vDal63{]} D. van Dalen, Extension problems in intuitionistic plane
projective geometry I, II, \emph{Indag. Math.} 25:349-383, 1963.

{[}vDal96{]} D. van Dalen, 'Outside' as a primitive notion in constructive
projective geometry, \emph{Geom. Dedicata} 60:107-111, 1996.  

{[}vPla95{]} J. von Plato, The axioms of constructive geometry, \emph{Ann.
Pure Appl. Logic} 76:169-200, 1995. 

{[}vPla98{]} J. von Plato, A constructive theory of ordered affine
geometry, \emph{Indag. Math.} \emph{(N.S.)} 9:549-562, 1998. 

{[}vPla10{]} J. von Plato, Combinatorial analysis of proofs in projective
and affine geometry, \emph{Ann. Pure Appl. Logic} 162:144-161, 2010.

{[}vSta47{]} K. G. C. von Staudt,\emph{ Geometrie der Lage, }Nürnberg,
1847. 

{[}VY10{]} O. Veblen and J. W. Young, \emph{Projective Geometry, Vol.
1,} Ginn, Boston, 1910. 

{[}Wei07{]} C. Weibel, Survey of non-Desarguesian planes, \emph{Notices
Amer. Math. Soc.} 54:1294\textendash 1303, 2007. 

{[}Y30{]} J. W. Young, \emph{Projective Geometry,} Math. Assoc. Amer.,
Open Court, Chicago, 1930. \\

2019 December 27 
\end{document}